\documentclass{amsart}

\usepackage{amscd,amsthm,xy,amssymb,nicefrac}
\xyoption{all}
\newdir{ >}{{}*!/-5pt/@{>}}

\newcommand{\com}{\mathfrak{K}}
\newcommand{\T}{\mathfrak{T}}
\newcommand{\A}{\mathbb{A}}

\newcommand{\B}{\mathcal{B}}
\newcommand{\I}{\mathcal{I}}
\newcommand{\J}{\mathcal{J}}

\newcommand{\F}{\mathcal{F}}

\newcommand{\M}{\mathcal{M}}
\newcommand{\cL}{\mathcal{L}}

\renewcommand{\inf}{\rm inf}
\newcommand{\alg}{\rm alg}
\newcommand{\dif}{\rm dif}
\newcommand{\nil}{\rm nil}
\renewcommand{\top}{\rm top}
\newcommand{\rel}{\rm rel}
\newcommand{\ran}{\rm ran}

\newcommand{\iso}{\overset{\cong}{\to}}
\newcommand{\weq}{\overset{\sim}{\to}}

\newcommand{\fib}{\twoheadrightarrow}

\newcommand{\triqui}{\triangleleft}

\newcommand{\locass}{\mathfrak{LocAlg}}
\newcommand{\ass}{\mathfrak{Ass}}
\newcommand{\ab}{\mathfrak{Ab}}
\newcommand{\hilb}{\mathfrak{Hilb}}

\newcommand{\lc}{{\rm top}}
\renewcommand{\mod}{{\rm Mod}}

\DeclareMathOperator{\tor}{Tor}

\DeclareMathOperator*{\colim}{colim}
\DeclareMathOperator*{\coli}{\colim}

\DeclareMathOperator{\LCV}{LC-\mathfrak{Vec}}
\DeclareMathOperator{\hofi}{hofiber}
\DeclareMathOperator*{\holi}{holim}
\newcommand{\DA}{\Delta^{\alg}}
\newcommand{\DD}{\Delta^{\dif}}
\newcommand{\hotimes}{\hat{\otimes}}

\newcommand{\circtimes}{\otimes_2}

\newcommand{\N}{\mathbb{N}}
\newcommand{\Z}{\mathbb{Z}}
\newcommand{\Q}{\mathbb{Q}}
\newcommand{\R}{\mathbb{R}}
\newcommand{\C}{\mathbb{C}}
\newcommand{\K}{\mathbb{K}}
\newcommand{\KH}{\mathbb{KH}}
\newcommand{\KD}{\mathbb{KD}}
\newcommand{\HN}{\mathbb{HN}}

\newcommand{\HP}{\mathbb{HP}}

\newcommand{\n}{\par\noindent}

\newcommand{\sn}{\par\smallskip\noindent}
\newcommand{\mn}{\par\medskip\noindent}

\theoremstyle{plain}
\newtheorem{thm}{Theorem}[subsection]
\newtheorem{coro}[thm]{Corollary}
\newtheorem{lem}[thm]{Lemma}
\newtheorem{propo}[thm]{Proposition}

\theoremstyle{definition}
\newtheorem{defi}[thm]{Definition}

\newtheorem{exa}[thm]{Example}

\theoremstyle{remark}
\newtheorem{nota}[thm]{Notations}
\newtheorem{rem}[thm]{Remark}

\begin{document}
\title[Algebraic and topological $K$-theory of locally convex algebras]{Comparison between algebraic and topological $K$-theory of locally convex algebras}

\author{Guillermo Corti\~nas}
\address{Guillermo Corti\~nas, Dep. Matem\'atica\\ Ciudad Universitaria Pab 1\\ (1428) Buenos Aires, Argentina\\ and  Dep. \'Algebra\\ Fac. de Ciencias\\
Prado de la Magdalena s/n\\ 47005 Valladolid, Spain.}
\email{gcorti@agt.uva.es} \urladdr{http://mate.dm.uba.ar/\~{}gcorti}

\author{Andreas Thom}
\address{Andreas Thom\\ Mathematisches Institut\\
Bunsenstr. 3-5\\ 37073 G\"ottingen, Germany.}
\email{thom@uni-math.gwdg.de}
\urladdr{http://www.uni-math.gwdg.de/thom/}
\begin{abstract}
This paper is concerned with the algebraic $K$-theory of locally convex $\C$-algebras stabilized by operator ideals,
and its comparison with topological $K$-theory. We show that if $L$ is locally convex and $\J$ a Fr\'echet operator ideal, then
all the different variants of topological $K$-theory agree on the completed projective tensor product $L\hotimes\J$, and that the obstruction for the comparison map
$K(L\hotimes\J)\to K^{\top}(L\hotimes\J)$ to be an isomorphism
is (absolute) algebraic cyclic homology. We prove the existence of an exact sequence (Theorem \ref{thm:main})
\[
\xymatrix{ K^{\top}_{1}(L\hotimes \J)\ar[r]&HC_{2n-1}(L\hotimes \J)\ar[r]&K_{2n}(L\hotimes \J)\ar[d] \\
K_{2n-1}(L \hotimes \J) \ar[u]& HC_{2n-2}(L \hotimes \J) \ar[l] & K_0^{\top}(L\hotimes \J). \ar[l]}
\]
We show that
cyclic homology vanishes in the case when $\J$ is the ideal of compact operators and $L$ is a
Fr\'echet algebra whose topology is generated by a countable family of sub-multiplicative seminorms and admits an approximate right
or left unit which is totally bounded with respect to that family (Theorem \ref{thm:karconj}).
This proves the generalized version of Karoubi's conjecture due to Mariusz Wodzicki and announced in
his paper {\it Algebraic $K$-theory and functional analysis}, First European Congress of Mathematics,
Vol. II (Paris, 1992), 485--496, Progr. Math., 120, Birkh\"auser, Basel, 1994.

We also consider stabilization with respect to a wider class of operator ideals, called sub-harmonic. Every Fr\'echet ideal
is sub-harmonic, but not conversely; for example the Schatten ideal $\cL_p$ is sub-harmonic for all $p>0$
but is Fr\'echet only if $p\ge 1$. We prove a variant of the exact sequence above which essentially says that if $A$ is a
$\C$-algebra and $\J$ is sub-harmonic, then the obstruction for the periodicity of $K_*(A\otimes_\C\J)$ is again cyclic homology
(Theorem \ref{thm:absolute_wodzicki}). This generalizes to all algebras a result of Wodzicki for $H$-unital algebras announced
in {\it loc.\,cit.}

The main technical tools we use are the diffeotopy invariance theorem of Cuntz and the second author (which we generalize in
Theorem \ref{thm:frech_hit}), and the excision theorem for infinitesimal $K$-theory, due to the first author.
\end{abstract}
\thanks{Corti\~nas' research was partly supported by grants PICT03-12330, UBACyT-X294, VA091A05, and
MTM00958. Thom's research was partly supported by the DFG (GK \textit{Gruppen und Geometrie} G\"ottingen).}

\maketitle

\tableofcontents
\section{Introduction}
This paper is about the comparison between algebraic and topological $K$-theory. This is a classical subject with numerous
applications, which has been considered by several authors, for different classes of topological algebras, and using a wide
variety of tools (see J.M. Rosenberg's excellent survey \cite{rosurvey}).

Here we are concerned with the comparison between algebraic and topological $K$-theory of not necessarily unital
locally convex $\C$-algebras. By a locally convex algebra we understand a complete locally convex vectorspace $L$ together with an
associative multiplication map $L\hotimes L\to L$; here $\hotimes$ is the projective tensor product of A. Grothendieck. We establish a six-term exact
sequence relating algebraic and topological $K$-theory with algebraic cyclic homology. We show (Theorem \ref{thm:exact}) that if $\J$ is a Fr\'echet operator
ideal and $L$ a locally convex algebra, then there is an exact sequence
\begin{equation}\label{intro:seqtop}
\xymatrix{ K^{\top}_{1}(L\hotimes \J)\ar[r]&HC_{2n-1}(L\hotimes \J)\ar[r]&K_{2n}(L\hotimes \J)\ar[d] \\
K_{2n-1}(L \hotimes \J) \ar[u]& HC_{2n-2}(L \hotimes \J) \ar[l] & K_0^{\top}(L\hotimes \J). \ar[l]}
\end{equation}
Here $K_*$ is algebraic $K$-theory. There are several possible definitions for topological $K$-theory of general locally convex algebras;
however we show (in Theorem \ref{thm:main}) that they all coincide for algebras of the form $L\hotimes\J$ as above. Thus for example in the sequence above we can define
$K^{\top}$ in terms of Cuntz' bivariant $K$-theory for locally convex algebras, see \cite{cw},
\[
K^{\top}_*(L\hotimes\J)=kk^{\top}_*(\C,L\hotimes\J).
\]
The algebraic cyclic homology groups are taken over $\Q$; this means that the (algebraic) tensor products appearing in the complex we use
for defining $HC$ (there are several quasi-isomorphic such complexes) must be taken over $\Q$. For example,
we have $HC_*(A)=H_*(C^\lambda(A),b)$, where $C^\lambda$ is Connes' complex, see \cite{lod},
\[
C_n^\lambda(A)=(A^{\otimes_\Q ^{n+1}})_{\Z/\langle n+1\rangle}.
\]
The meaning of the sequence \eqref{intro:seqtop} is clear; it says that, for locally convex algebras stabilized by a Fr\'echet
operator ideal, algebraic cyclic homology measures the obstruction for the comparison map $K_*\to K_*^{\top}$ to be an isomorphism.
As an immediate application of this, and of the fact that, by definition, cyclic homology vanishes in negative degrees,
we get
\begin{equation}\label{intro:kneg}
K_n(L\hotimes\J)=K^{\top}_n(L\hotimes \J) \qquad (n\le 0).
\end{equation}
The particular case of \eqref{intro:kneg} when $\J=\cL_p$, and $p> 1$ (or, more generally, when $\J$ is harmonic) was proved
in \cite[Thm. 6.2.1]{ct}.
\sn
It is also clear from \eqref{intro:seqtop} that the vanishing of cyclic homology in {\it all} degrees is equivalent to the
isomorphism between algebraic and topological $K$-theory. For example, we show that if $\J=\cL_\infty$ is the ideal of all
compact operators and $L$ is a unital Banach algebra then $HC_*(L\hotimes\J)=0$, whence
\begin{equation}\label{intro:karconj}
K_*(L\hotimes\J)=K_*^{\top}(L\hotimes\J).
\end{equation}
This establishes Karoubi's conjecture (as stated in \cite{karcomp}). In fact we show \eqref{intro:karconj} holds more generally when $L$ is a
Fr\'echet algebra whose topology is generated by a countable family of sub-multiplicative seminorms and admits an approximate right
or left unit
which is totally bounded with respect to that family. We point out that Mariusz Wodzicki is credited with the solution of Karoubi's
conjecture, both the original one and the generalization just mentioned, as well as with other
results proved in this paper. He has lectured on these results in several places, including Heidelberg and Paris, giving
full details of his proofs. However, although some of these results have been announced in \cite{wod}, his proofs have not been
published in print except in some particular cases, see \cite{hus}. Our proofs use some of the published results of Wodzicki,
as well as other results which are independent of his work. For example, most of our proofs rely heavily on
the diffeotopy invariance theorem from \cite{ct} --which we generalize in \ref{thm:frech_hit}-- and the excision theorem for infinitesimal $K$-theory from \cite{kabi}, none of which were
available at the time when Wodzicki's pioneering work \cite{wod} appeared.
\sn
Another result announced in \cite[Thm. 5]{wod} is the existence of a $6$-term exact sequence
\begin{equation}\label{intro:rel_seq}
\xymatrix{ K_{-1} (A\otimes_\C \J)\ar[r]&HC_{2n-1}(A\otimes_\C\B:A\otimes_\C \J)\ar[r]&K_{2n}(A\otimes_\C\B:A\otimes_\C \J)\ar[d] \\
K_{2n-1}(A\otimes_\C\B:A\otimes_\C \J) \ar[u] & HC_{2n-2}(A\otimes_\C\B:A\otimes_\C \J) \ar[l] & K_{0}(A\otimes_\C \J). \ar[l]}
\end{equation}
Here $A$ is an $H$-unital $\C$-algebra, $\otimes_\C$ is the algebraic tensor product, $\B$ is the algebra of bounded operators
in an infinite dimensional, separable Hilbert space, and $\J$ is what in this paper we call a sub-harmonic
operator ideal (see \ref{defi:subhar} for a definition). We prove \eqref{intro:rel_seq} for all algebras $A$ and for
all sub-harmonic operator ideals $\J$. Thus we generalize Wodzicki's sequence from $H$-unital algebras to all algebras. Furthermore, we also prove a
variant of \eqref{intro:rel_seq}, which is still valid under the same hypothesis, and which involves absolute, rather than relative
$K$-theory and cyclic homology. We show that there is an exact sequence
\begin{equation}\label{intro:abs_seq}
\xymatrix{ K_{-1} (A\otimes_\C \J)\ar[r]&HC_{2n-1}(A\otimes_\C \J)\ar[r]&K_{2n}(A\otimes_\C \J)\ar[d] \\
K_{2n-1}(A \otimes_\C \J) \ar[u] & HC_{2n-2}(A \otimes_\C \J) \ar[l] & K_{0}(A\otimes_\C \J). \ar[l]}
\end{equation}
Examples of sub-harmonic ideals include all Fr\'echet ideals (\ref{propo:subharm}) as well as some ideals, such as the Schatten
ideals $\cL_p$ with $0<p<1$, which are not even locally convex. In the particular case when $A=\C$ both \eqref{intro:rel_seq} and
\eqref{intro:abs_seq} simplify. Indeed we show in \ref{thm:compk0k1} that, as stated without proof in \cite[Prop. on p. 491]{wod},
we have
\begin{equation}\label{intro:compk0k1}
K_{-1}(\J)=0,\ \ K_0(\J)=\Z
\end{equation}
for any proper operator ideal $\J$. Thus \eqref{intro:abs_seq} becomes
\begin{equation}\label{intro:absac}
0\to HC_{2n-1}(\J)\to K_{2n}(\J)\to \Z\overset{\alpha_n}\to HC_{2n-2}(\J)\to K_{2n-1}(\J)\to 0.
\end{equation}
We show moreover that if $\I\subset\cL_p$ $(p\ge 1)$ then $\alpha_n$ is injective for $n\ge (p+1)/2$. As an application of this
we obtain in \ref{subsec:ck} a new description of the multiplicative character of a $p$-summable Fredholm module defined by
A. Connes and M. Karoubi in \cite{ck}.
\sn
One can also combine \eqref{intro:compk0k1} with \eqref{intro:rel_seq} to obtain a sequence similar to \eqref{intro:absac},
but involving relative, instead of absolute $K$-theory and cyclic homology:
\begin{equation}\label{intro:relac}
0\to HC_{2n-1}(\J)\to K_{2n}(\B:\J)\to \Z\to HC_{2n-2}(\B:\J)\to K_{2n-1}(\J)\to 0.
\end{equation}
This is the sequence announced in \cite[Thm. 6]{wod}; Dale Husem\"oller took the work to write down Wodzicki's proof in \cite{hus}.
\sn
As indicated above we obtain a generalization of the diffeotopy invariance theorem proved by J. Cuntz and the second author in \cite{ct}. The latter implies that if
$E$ is a functor from the category $\locass$ of locally convex algebras to abelian groups which is split exact and stable
with respect to $2\times 2$-matrices, and $\J$ is a harmonic ideal, then
\begin{equation}\label{intro:frech_hit}
L\mapsto E(L\hotimes\J)
\end{equation}
is diffeotopy invariant, i.e. sends the two evaluation maps $A \hotimes C^{\infty}(\Delta^1) \to A$ to the same morphism. Recall that a harmonic ideal is a Banach operator ideal $\J\subset\B$ with continuous inclusion, which
contains a compact operator whose sequence of singular values is the harmonic sequence $(\nicefrac1n)_n$, and which is multiplicative
(see \ref{subsec:opideals} for the relevant definitions). We prove that, under an extra hypothesis on $E$, the functor
\eqref{intro:frech_hit}
is diffeotopy invariant for any Fr\'echet ideal $\J$. The extra hypothesis essentially says that $E$ sends certain homomorphisms whose kernel and cokernel are both square-zero algebras into isomorphisms (see \ref{defi:nilinv} for
a precise definition). For example the functors $K_n$ for $n\le 0$ as well as the infinitesimal and (polynomial) homotopy $K$-theory groups
$K^{\inf}_n$ and $KH_n$ for $n\in\Z$, all satisfy these hypothesis. Hence if $A$ is any $\C$-algebra, and $\J$ any Fr\'echet
ideal, then the functors
\begin{equation}\label{intro:maini}
KH_*(A\otimes(?\hotimes\J))\text{ and } K^{\inf}_*(A\otimes(?\hotimes\J)),
\end{equation}
are diffeotopy invariant. Using this we show, for example, that
\begin{equation}\label{intro:reg}
A\otimes(L\hotimes\J)\text{ is $K_0$-regular.}
\end{equation}
In particular
\[
K_n(A\otimes(L\hotimes\J))=KH_n(A\otimes(L\hotimes\J))\ \ (n\le 0).
\]
Further, we prove that topological and homotopy
$K$-theory agree on stable locally convex algebras. We have
\begin{equation}\label{intro:kh=ktop}
KH_*(L\hotimes\J)=K^{\top}_*(L\hotimes\J)
\end{equation}
for every Fr\'echet ideal $\J$. In particular both \eqref{intro:reg} and \eqref{intro:kh=ktop} hold when $\J=\cL_p$, $p\ge 1$.
\sn
We also correct a mistake in \cite[Thm 7.1]{ct}, which states that if $\J$ is a harmonic operator ideal,
then the map
\begin{equation}\label{intro:det}
\xymatrix{K_1(\J)\ar[r]& K_1(\J/\J^2)\ar[r]^(0.6){det}& \J/\J^2}
\end{equation}
is an isomorphism. Actually the proof given in {\it loc.\ cit.} only works under the extra assumption that also $\J^2$
be harmonic (see \ref{rem:k1harmonic} for details). In fact, we show in \ref{rem:pearcy-topping}, using results from
\cite{dykwodz}, that the map \eqref{intro:det} is not injective for $\J=\cL_2$.
\mn
The rest of this paper is organized as follows. In Section \ref{sec:hc} we recall some basic facts about the $K$-theory and cyclic homology
of $\Q$-algebras, their different variants, and the Chern characters between them. Some of these are displayed in the following map of
exact sequences ($n\in\Z$)
\begin{equation}\label{intro:fundechar}
\xymatrix{ KH_{n+1}A\ar[r]\ar[d]_{ch_{n+1}}& K^{\nil}_nA\ar[r]\ar[d]_{\nu_n}&K_n A\ar[r]\ar[d]_{c_n}& KH_nA
\ar[d]_{ch_n}\ar[r]&\tau K^{\nil}_{n-1}A\ar[d]_{\nu_{n-1}}\\
               HP_{n+1}A\ar[r]_S&HC_{n-1}A\ar[r]_B&HN_nA\ar[r]_I&HP_nA\ar[r]_S&HC_{n-2}A.}
\end{equation}
The basic or primary character is the Jones-Goodwillie map $c_*$; it takes values in negative cyclic homology. Its (polynomial)
homotopy invariant version, $ch_*$, goes from homotopy $K$-theory to periodic cyclic homology. The map $\nu_*$ is the secondary character; it goes from
$\nil$-$K$-theory to cyclic homology. We give a reasonably self-contained exposition of how these maps are constructed.
The construction we give combines methods of \cite{crelle} and \cite{wenil}.
We do not claim much originality. Actually a diagram similar to \eqref{intro:fundechar}, involving Karoubi-Villamayor $K$-theory $KV$ instead of $KH$ (which had not yet been invented by C. Weibel), appeared
in \cite{wenil} (see also \cite{karmult}). For $K_0$-regular algebras and $n\ge 1$, the latter diagram is equivalent to
\eqref{intro:fundechar} (see Subsection \ref{subsec:compakw}).
The vertical maps in \eqref{intro:fundechar} are
induced by maps of spectra after taking homotopy groups; the homotopy groups of their fibers form a long exact sequence
\begin{equation}\label{intro:infchar}
\xymatrix{ K^{\inf}H_{n+1}A\ar[r]& K^{\inf,\nil}_nA\ar[r]&K^{\inf}_n A\ar[r]& K^{\inf}H_nA.}
\end{equation}
Here $K^{\inf}$ is infinitesimal $K$-theory and $K^{\inf}H$ its polynomial homotopy variant; the groups $K_*^{\inf,\nil}$ are the relative groups.
In the last subsection of
this section, we give an alternative proof, for the case of $\Q$-algebras, of Weibel's theorem that $KH$ satisfies excision.
The proof is based on the diagram \eqref{intro:fundechar} and on the excision properties of $K^{\inf}$ (\cite{kabi})
and $HP$ (\cite{cq}).
\sn
Section \ref{sec:kinfreg} is concerned with polynomial homotopy invariance
for infinitesimal $K$-theory, $K^{\inf}$. Recall that a if $F$ is a functor from algebras to abelian groups, then an algebra
$R$ is called $F$-regular if for the polynomial ring, the map $F(R)\to F(R[t_1,\dots,t_m])$ is an isomorphism for each $m\ge 1$.
T. Vorst has shown in \cite{vorst} that $K_n$-regularity implies $K_{n-1}$-regularity. We prove in Proposition \ref{propo:vokinf}
that the same is true of $K^{\inf}$; a $K_n^{\inf}$-regular $\Q$-algebra is $K^{\inf}_{n-1}$-regular. Using this together with
\eqref{intro:infchar} and \eqref{intro:fundechar}, we prove in \ref{propo:nuiso} that if $A$ is $K^{\inf}_n$-regular, then
the secondary character $\nu_m$ is an isomorphism for $m\le n$. In particular if $A$ is $K^{\inf}$-regular (i.e. $K_n^{\inf}$-regular
for all $n$) then
\begin{equation}\label{intro:nuiso}
\nu_*:K_*(A)\iso HC_{*-1}(A)\qquad (A \ \ K^{\inf}\text{-regular}).
\end{equation}
In Subsection \ref{subsec:regkink} we compare $K^{\inf}$-regularity with $K$-regularity; for example we show that $K^{\inf}_0$-regular
algebras are $K_0$-regular. Some examples of $K^{\inf}$-regular algebras are given in Subsection \ref{subsec:exaregkinf}; they include
nilpotent algebras and also Wagoner's infinite sum algebras, such as the algebra $\B$ of bounded operators in an infinite dimensional
Hilbert space (see Example \ref{exa:nilpo} and Lemma \ref{lem:infisum}). On the other hand we show that no nonzero
unital commutative algebra can be $K^{\inf}$-regular (Proposition \ref{propo:notcom}).
\sn
In Section \ref{sec:locass} we construct, for locally convex algebras $A$, a diagram similar to \eqref{intro:fundechar}, with
diffeotopy invariant $K$-theory $KD$ substituted for $KH$, and the continous or topological versions of cyclic homology and
its variants substituted for their algebraic counterparts:
\begin{equation}\label{intro:diffechar}
\xymatrix{KD_{n+1}A\ar[d]_{ch^{\dif}_{n+1}}\ar[r]&K^{\rel}_nA\ar[d]_{ch^{\rel}_n}\ar[r]&K_nA\ar[d]_{\hat{c}_n}\ar[r]&KD_nA\ar[d]_{ch^{\dif}_n}\\           HP^{\top}_{n+1}A\ar[r]_S&HC^{\top}_{n-1}A\ar[r]_B&HN^{\top}_n A\ar[r]_I&HP_n^{\top}A}
\end{equation}
The construction is analogous to that of \eqref{intro:fundechar}. Again, we do not claim much originality, since a version of this
for the connected, Karoubi-Villamayor variant of $KD$ was constructed by Karoubi in \cite{karmult} for Fr\'echet algebras.
When $A$ is well-behaved, so that $K_n(A)$ is diffeotopy invariant for $n\le 0$, Karoubi's diagram coincides with
\eqref{intro:diffechar} for $n\ge 1$ (see Subsection \ref{subsec:compaktop}). Moreover, there is a map of diagrams going from \eqref{intro:fundechar} to \eqref{intro:diffechar}.
Thus whenever the comparison map $KH_*(A)\to KD_*(A)$, is an isomorphism, the top rows identify, and the topological characters
factor through the algebraic ones. Also in this section we prove (by a similar method as that sketched above for $KH$) that $KD$ satisfies excision (\ref{propo:kdexc}).
\sn
In Section \ref{sec:opid_stab} we fix our notations for operator ideals, and prove some basic properties on them. For us an operator
ideal is a functor defined on the category $\hilb$ of infinite dimensional, separable $\C$-Hilbert spaces and isometries, which associates
to each $H\in\hilb$ a proper ideal $\J(H)$ of the algebra $\B(H)$ of bounded linear operators. We show in Proposition
\ref{propo:subharm} that if $\J$ is a complete locally convex operator ideal such that the multiplication
map $\B\times\J\times\B\to\J$ is jointly continous, then $\J\supset\cL_1$, and moreover, for the completed tensor product $\circtimes$
of Hilbert spaces, the map $\boxtimes:B(H)\otimes \B(H)\to \B(H\circtimes H)$ sends $\cL_1\otimes\J$ into $\J$. In symbols
\[
\cL_1(H)\boxtimes\J(H)\subset\J(H\circtimes H).
\]
Because of this, we say that $\cL_1$ multiplies $\J$. For example every Fr\'echet (operator) ideal $\J$ with continuous inclusion
$\J\subset\B$ satisfies the hypotesis of Proposition \ref{propo:subharm} (see Subsection \ref{subsec:opideals}); in particular it contains and is multiplied by $\cL_1$.
In Subsection \ref{subsec:perio} we show (Proposition \ref{propo:periodicity}) that if $\J$ is a Fr\'echet ideal, $X$ any excisive, diffeotopy invariant, $M_2$-stable
functor from locally convex algebras to spectra, then $X(?\hotimes\J)$ is $2$-periodic; thus
$X(?\hotimes\J)\weq \Omega^2(X(?\hotimes\J))$ for every locally convex algebra $L$.
As an application, we show that
$KD_*(A\hotimes\J)=KD_{*+2}(A\hotimes\J)$ for any Fr\'echet operator ideal $\J$. In Subsection \ref{subsec:kk} we recall --from
\cite{cw},\cite{ct}, and \cite{wt}-- some fundamental
facts  about J. Cuntz' bivariant $K$-theory for locally convex algebras, $kk^{\top}$, and about its algebraic counterpart, $kk$.
\sn
The main results of the paper are contained in Section \ref{sec:main}. The first of these is the generalized version
\eqref{intro:frech_hit} of the Cuntz-Thom diffeotopy invariance theorem of \cite{ct}, which is valid for stabilization with Fr\'echet ideals
(Theorem \ref{thm:frech_hit}). Next we prove Theorem \ref{thm:main}, which establishes several properties
of $K$-theory stablilized by Fr\'echet ideals. For every
$\C$-algebra $A$, locally convex algebra $L$ and Fr\'echet ideal $\J$, we show that the functors \eqref{intro:maini} are diffeotopy invariant,
that $A\otimes_\C(L\hotimes\J)$ is $K^{\inf}$-regular, and that the different variants of topological $K$-theory agree after stabilizing; for
example:
\begin{equation}\label{intro:allagree}
KD_n(L\hotimes \J)=kk^{\top}_n(\C,L\hotimes \J)=KH_n(L\hotimes \J) \quad \forall n\in\Z.
\end{equation}
Note that \eqref{intro:kneg} follows from this and the cited $K^{\inf}$-regularity of $L\hotimes\J$. Indeed, as mentioned above,
$K^{\inf}$-regularity implies $K_0$-regularity; on the other hand the latter implies the agreement between $K_n$ and $KH_n$ for nonpositive $n$.
Another result proved in this section is the exact sequence \eqref{intro:seqtop} (Theorem \ref{thm:exact}).
The proof uses the $K^{\inf}$-regularity of $L\hotimes\J$ proved in Theorem \ref{thm:main},
together with \eqref{intro:fundechar}, \eqref{intro:nuiso}, and \eqref{intro:allagree}. The computation of $K_0$ and $K_{-1}$ for
arbitrary operator ideals \eqref{intro:compk0k1} is the subject of Theorem \ref{thm:compk0k1}.
In Subsection \ref{subsec:main_alg} we introduce the notion of sub-harmonic ideal. An operator ideal $\J$ is sub-harmonic if there
exists a $p>0$ such that $\cL_p$ multiplies the root completion
\[
\J_{\infty} = \bigcup_{n\geq 1} \J^{\nicefrac1n}.
\]
We show (Theorem \ref{thm:main_alg}) that if $\J$ is sub-harmonic and $A$ is any $\C$-algebra, then $A\otimes_\C\J$ is
$K^{\inf}$-regular, and $KH_*(A\otimes_\C\J)$ is $2$-periodic:
\[
KH_n(A \otimes_\C\J)=\left\{\begin{matrix} K_0(A\otimes_\C\J)& \text{ if }n\text{ is even.}\\
                                           K_{-1}(A\otimes_\C\J)&\text{ if }n\text{ is odd.}\end{matrix}\right.
\]
\sn
Section \ref{sec:appli} is devoted to some applications of the results of the previous section. These include
the exact sequences \eqref{intro:rel_seq} and \eqref{intro:abs_seq} (Theorem \ref{thm:absolute_wodzicki}) as well as \eqref{intro:absac}
(Proposition \ref{propo:kopi}), a computation of the $K$-theory
of a generalized version of the representation algebras considered by Connes and Karoubi in their construction of the multiplicative
character of a Fredholm module (Proposition \ref{propo:mji}), and other results on the lower $K$-theory
and cyclic homology of operator ideals, including
the correction of \cite[Thm. 7.1]{ct} cited above for the case when $\J$ is harmonic (\ref{rem:k1harmonic}-\ref{rem:pearcy-topping}).
\sn
Section \ref{sec:karconj} is devoted to Karoubi's conjecture. After some basic facts on $H$-unitality and Hochschild homology, which
are the subject of Subsection \ref{subsec:homprep}, we give a proof of Wodzicki's result (stated without proof in \cite[Thm. 4]{wod})
that, for an operator ideal $\J$, the condition that $\J=\J^2$ is equivalent to the $H$-unitality of
$\J$ as a $\Q$-algebra (Theorem \ref{thm:jhunit}). In Theorem \ref{thm:alg_carcoj} we prove that is sub-harmonic and $\J=[\J,\J]$ and $A$ is
an $H$-unital $\C$-algebra, then $HC_*(A\otimes_\C\J)=0$, $A\otimes_\C\J$ is $K$-regular, and $K_*(A\otimes_\C\J)$ is $2$-periodic.
The particular case of this periodicity result when $\J=\cL_\infty$ is due to Wodzicki and stated without proof in \cite[Thm. 2 (c)]{wod}.
In Theorem \ref{thm:karconj} we prove Karoubi's conjecture for Fr\'echet algebras whose topology is generated by a countable family of sub-multiplicative seminorms and admits an approximate
right or left unit
which is totally bounded with respect to that family, another result of Wodzicki's, stated without proof in \cite[Thm. 2]{wod}.
The proof uses the homology vanishing results of previous subsections, together with sequence \eqref{intro:seqtop}
(Theorem \ref{thm:exact}) and Wodzicki's theorem (proved in \cite{wodex} and
recalled in \ref{thm:wodzh} below)
that Fr\'echet algebras satisfying the extra hypothesis of Theorem \ref{thm:karconj} are $H$-unital. An alternative proof,
using only Wodzicki's theorem \ref{thm:wodzh} and the homotopy invariance theorem of \cite{hig},
is given in Remark \ref{rem:kartriv}. However we point out --also in \ref{rem:kartriv}-- that the proof given in
\ref{thm:karconj} gives hope that it may be adapted to prove an appropriate variant of Karoubi's conjecture for algebras
stabilized with other ideals $\J$ such that $\J=[\J,\J]$ instead of $\cL_\infty$, such as the root completion $(\cL_1)_\infty$
of the ideal of trace class operators $\cL_1$.
\mn
\n{\bf Acknowledgements.} This paper grew out of the individual work of the first author on the sequence \eqref{intro:seqtop}, of
the second on generalizing the diffeotopy invariance theorem \eqref{intro:frech_hit}, and of our joint efforts in trying to prove the
results announced in \cite{wod}. We finally decided to put all together in one article, which showed to be a good idea, since
in writing it up we obtained stronger results than we initially had separately. Part of this research was carried out during
a visit of the second author to the University of Valladolid; he is thankful to this institution for its hospitality.
\sn
We both wish to thank Mariusz Wodzicki, for his work, and in particular his paper \cite{wod}, was largely responsible for motivating us in this research.

\section{Cyclic homology}\label{sec:hc}

\subsection{Preliminaries}

All rings considered in this paper shall be (not necessarily unital) algebras
over the field $\mathbb{Q}$ of rational numbers. We write $C$, $CC$, $CN$ and $CP$ for the Hochschild, cyclic, negative
cyclic, and periodic cyclic complexes; thus $HH_*A=H_*C(A)$,
$HC_*A=H_*CC(A)$, $HN_*A=H_*CN(A)$ and
$HP_*A=H_*CP(A)$ are {\it Hochschild, cyclic,
negative cyclic and periodic cyclic homology} groups of $A$. There is
a commutative diagram with exact rows
\[
\xymatrix{0\ar[r]&CN(A)\ar[d]\ar[r]^I&CP(A)\ar[d]^S\ar[r]^S&CC(A)[-2]\ar[r]\ar[d]^1&0\\
          0\ar[r]&C(A)\ar[r]_I& CC(A)\ar[r]_S&CC(A)[-2]\ar[r]&0}
\]
which gives rise to a map of exact sequences
\begin{equation}\label{sbis}
\xymatrix{HP_{n+1}A\ar[d]\ar[r]^S&HC_{n-1}A\ar[r]^B\ar[d]^1&HN_nA\ar[r]^I\ar[d]&HP_nA\ar[d]\ar[r]^S&HC_{n-2}A\ar[d]^1\\
HC_{n+1}A\ar[r]_S&HC_{n-1}A\ar[r]_B&HH_nA\ar[r]_I&HC_nA\ar[r]_S&HC_{n-2}.}
\end{equation}
The top and bottom row in the diagram above are Connes' {\it $SBI$-sequences}.
There is an obvious abuse of notation in that several distinct arrows are given the same name, but this is the standard notation.

The main properties of $HP$ which we shall use are that it maps short exact sequences of algebras to long exact sequences ({\it excision}), that it vanishes
on nilpotent algebras ({\it nilinvariance}) and that it maps the inclusion
$A\subset A[t]$ to an isomorphism ({\it polynomial homotopy invariance}).
The first of these results is due to Cuntz and Quillen \cite{cq} and the other
two to Goodwillie \cite{goo}, \cite{goo1}.

Each of the homologies considered above extend to
simplicial algebras in the standard way, by applying the corresponding
functorial complex degreewise and then taking the total complex. Next we recall a few results about
homology of simplicial algebras which we shall need. Although these are well-known (e.g. see \cite{gw}),
we include short proofs of each of them.

\begin{lem}\label{pinulhcnul}
Let $A$ be a simplicial algebra, $\pi_*A$ its homotopy groups. Assume
$\pi_nA=0$ for all $n$. Then $HH_*A=HC_*A=0$ and $HN_*A=HP_*A$.
\end{lem}
\begin{proof}
Let $p,q\ge 0$; recall that $C_qA_p=(A_p\oplus \Q)\otimes A_p^{\otimes q}$. By K\"unneth's theorem,
$C_q(A_*)$ is contractible for $q \geq 0$.
The assertion for Hochschild homology is immediate from this; the remaining assertions follow using the $SBI$-sequences.
\end{proof}

Let $A$ be an algebra. Consider the simplicial algebra
\[
\DA A:[n]\mapsto A\otimes\Q[t_0,\dots,t_n]/\langle1-(t_0+\dots+t_n)\rangle.
\]
We identify elements of $\DA_nA$ with polynomial functions on
the algebraic $n$-simplex $\{(x_0,\dots,x_n)\in \Q^{n+1}:\sum x_i=1\}$ with
values in $A$. Face and degeneracy maps are given by
\begin{gather}\label{formucaradege}
d_i(f)(t_0,\dots,t_{n-1})=f(t_0,\dots,t_{i-1},0,t_i,\dots,t_{n})\\
s_j(f)(t_0,\dots,t_{n+1})=f(t_0,\dots,t_{i-1},t_i+t_{i+1},\dots,t_{n+1}).\nonumber
\end{gather}
Here $f\in \DA_nA$, $0\le i\le n$, and $0\le j\le n-1$.

In the course of the proof of the lemma below and elsewhere, we consider the {\it Moore complex} $\M V$ of a simplicial abelian group $V$.
Recall that $\M _nV:=\bigcap_{i=0}^{n-1}\ker d_i$ and that the boundary map is
$d_n:\M_n V\to \M_{n-1}V$. The main property of $\M V$ is that its homology groups are the homotopy groups of $V$.
\begin{lem}
$\pi_*\DA A=0$.
\end{lem}
\begin{proof} By the Dold-Kan correspondence, the homotopy groups of
 $\DA A$ are the homology groups of the Moore complex $\M_*:=\M_*(\DA A)$.
Note $\M_n\subset \DA_nA$ consists of those
polynomials $f$ which vanish on all faces of the $n$-simplex except perhaps on the
$n$th face. Define inductively a map $h_n:\M_n\to \M_{n+1}$ as follows
\begin{equation}\label{hcontra}
h_n(f)=(1-t_{n+1})(f-h_{n-1}d_nf).
\end{equation}
One checks that $1=d_{n+1}h_n+h_{n-1}d_n$,
whence $h_*$ is a contracting homotopy for $\M_*$.
\end{proof}
\begin{coro}\label{bafa}
The maps $HN_*(\DA A)\to HP_*(\DA A)\leftarrow HP_*(A)$ are
isomorphisms.
\end{coro}
\begin{proof}
That the first map is an isomorphism follows from the lemmas above
and Connes' $SBI$-sequences \eqref{sbis}. That also the second map is an isomorphism
is a consequence of the polynomial homotopy invariance of $HP$.
\end{proof}

\subsection{Primary and secondary Chern characters}\label{pschalg}
Let $A$ be a unital algebra. We write $BG$ for the nerve of the group
$G$; thus for us $BG$ is a pointed simplicial set. A functorial model
for the plus construction of the general linear group is the Bousfield-Kan
$\Z$-completion
$$
K(A):=\Z_{\infty}BGL(A).
$$
There is a nonconnective spectrum $\K A$ of which the $n$-th space
is
\begin{equation}\label{nspace}
{}_n\K A:=\Omega K(\Sigma^{n+1}A)
\end{equation}
where $\Sigma$ is Karoubi's suspension functor (\cite{gersten}).
The basic or {\it primary} character from $K$-theory to cyclic homology is the
{\it Jones-Goodwillie Chern character}
\begin{equation}\label{chee}
c_n:K_n(A)\to HN_n(A)
\end{equation}
defined for all $n\in \Z$. For $n\ge 1$ it is
induced by a map of spaces
\begin{equation}\label{chespace}
c:K (A)\to SCN_{\ge 1}(A).
\end{equation}
Here $SCN_{\ge 1}$ is the simplicial abelian group the Dold-Kan
correspondence associates to the truncation of $CN$ which is $0$
in degree $0$, the kernel of the boundary operator in degree one,
and $CN_n$ in degrees $n\ge 2$. If $n\le 0$, then applying \eqref{chespace} to
$\Sigma^{n+1}A$, using the equivalence
\begin{equation}\label{listin}
S(CN_{\ge 1-(n+r)}(A)[-(n+r)])\weq SCN_{\ge 1}(\Sigma^{n+r}A)
\end{equation}
(see \cite[\S4]{kabi}) and taking homotopy groups, one gets the map
\[
c_n:K_nA=\pi_0({}_n\K A)\to\pi_0 SCN(\Sigma^{n+1}A)=HN_nA \ \ (n\le 0).
\]
\mn
The {\it infinitesimal $K$-theory spectrum} of $A\in\ass_1$  is
the fibrant spectrum
\begin{equation}\label{finfx}
\K^{\inf}A:=\holi_{A_{\inf}}\K.
\end{equation}
Here $A_{\inf}$ is the category of all surjections $B\fib A\in\ass_1$
with nilpotent kernel; $\K$ is viewed as a functor on $A_{\inf}$ by
$(B\to A)\mapsto \K(B)$. Although the category $A_{\inf}$ is large,
it is proved in \cite[5.1]{crelle} that it has a left cofinal small
subcategory in the sense of \cite[Ch IX\S9]{bk}, which by
\cite[2.2.1]{crelle} can be chosen to depend functorially on $A$.
The homotopy limit is taken over this small subcategory.

There is a natural map $\K^{\inf}A\to \K A$;
we write $\tau \K A$ for the delooping of its homotopy fiber.
We have a homotopy fibration
\begin{equation}\label{basicfib}
\xymatrix{\K^{\inf}A\ar[r]&\K A\ar[r]^{c^\tau}& \tau \K A.}
\end{equation}
We call $c^\tau$ the {\it tautological character}; it is closely related to the Jones-Goodwillie
character. Indeed it is shown in \cite[\S2]{kabi} that there is a natural isomorphism $HN_n(A)\iso\tau K_n(A)$
($n\in \Z$) which makes the following diagram commute
\begin{equation}\label{agripa}
\xymatrix{K_n(A)\ar[r]^{c^\tau_n}\ar[dr]^{c_n} &\tau K_n(A)\\
           &HN_n(A).\ar[u]^{\wr}}
\end{equation}
Moreover, we have the following.

\begin{lem}\label{hntauk}
Let ${}_n\tau \K A$ be the $n$th space of the $\tau \K A$ spectrum ($n\ge 0$).
Then there is a weak equivalence
$\Omega SCN_{\ge 1}(\Sigma^{n+1}A)\weq {}_n\tau\K A$ such that the following
diagram is homotopy commutative
\[
\xymatrix{{}_n\K A\ar[r]\ar[d]&{}_n\tau \K A\\
\Omega SCN_{\ge 1}(\Sigma^{n+1}A)\ar[ur]&}
\]
\end{lem}
\begin{proof}
In the proof of \eqref{agripa} given in \cite{kabi} it is shown that for each $n\ge 0$ there is a homotopy
commutative diagram
$$
\xymatrix{{}_{n-1}\tau \K(A)\ar[r]\ar[d]^{\wr}&\holi_{A_{\inf}}{}_{n}\K
\ar[r]\ar[d]&{}_n\K(A)\ar[d]\\
F \ar[r]&\holi_{A_{\inf}}\Omega SCN_{\ge 1}\Sigma^{n+1}\ar[r] &\Omega SCN_{\ge 1}(\Sigma^{n+1}A)}
$$
where the vertical map on the left is a weak equivalence and the middle space in the bottom line is contractible.
The lemma is immediate from this.
\end{proof}
Recall Weibel's  {\it (polynomial) homotopy} $K$-theory spectrum
$\KH(A)$ is defined as the total fibrant spectrum $\K(\DA A)$ associated to the simplicial
spectrum $[n]\mapsto \K(\DA_n A)$ (\cite{wh}). There is a natural map $\K A\to \KH A$; we write $\K^{\nil}A$ for its homotopy
fiber. We have a homotopy commutative diagram with homotopy fibration rows and columns
\begin{equation}\label{fundechar}
\xymatrix{\K^{\inf,\nil}A\ar[d]\ar[r]&\K^{\inf}A\ar[r]\ar[d]&\K^{\inf}(\DA A)\ar[d]\\
           \K^{\nil}A\ar[r]\ar[d]_\nu & \K A\ar[d]_{c^\tau}\ar[r]&\KH A\ar[d]_{ch^\tau}\\
           \tau\K^{\nil}A\ar[r]&\tau\K A\ar[r]&\tau\K(\DA A).}
\end{equation}
The middle column is \eqref{basicfib}; that on the right is \eqref{basicfib} applied to $\DA A$; the horizontal map of
homotopy fibrations from middle to right is induced by the inclusion $A\to \DA A$, and its fiber is the column on the left.

\begin{lem}\label{htodotautodo}
There are natural isomorphisms $$HC_{n-1}A\iso \pi_n\tau\K^{\nil}A, \quad \mbox{and}  \quad HP_nA\iso\pi_n\tau\K(\DA A), \quad \forall n\in\Z,$$ which,
together with that of \eqref{agripa}, define an isomorphism of long exact sequences
\[
\xymatrix{HP_{n+1}A\ar[r]_S\ar[d]_\wr&HC_{n-1}A\ar[r]_B\ar[d]_\wr& HN_nA\ar[r]_I\ar[d]_\wr&HP_n(\DA A)
\ar[d]_\wr\ar[r]_S&HC_{n-2}A\ar[d]_\wr\\
               \tau K_{n+1}(\DA A)\ar[r]&\tau K^{\nil}_nA\ar[r]&\tau K_nA\ar[r]&\tau K_nA\ar[r]&\tau K^{\nil}_{n-1}A}
\]
where the top line is the $SBI$-sequence of \eqref{sbis}.
\end{lem}
\begin{proof}
Let $\iota:A\to \DA A$ be the natural map. By \ref{hntauk},\eqref{listin} and \ref{bafa},
there is a homotopy commutative diagram in which all vertical arrows are weak equivalences
\[
\xymatrix{{}_n\tau\K(A)\ar[r]_\iota & {}_n\tau\K(\DA A)\\
           \Omega SCN_{\ge 1}(\Sigma^{n+1}A)\ar[u]^\wr\ar[r]_\iota &\Omega SCN_{\ge 1}(\Sigma^{n+1}\DA A)\ar[u]^\wr\\
           S(CN_{\ge -n}(A)[-n])\ar[r]_\iota\ar[u]_\wr&S(CN_{\ge -n}(\DA A)[-n])\ar[u]^\wr\ar[d]^I\\
           S(CN_{\ge -n}(A)[-n])\ar[r]_{I\iota}\ar[u]^1&S(CP_{\ge -n}(\DA A)[-n])\\
           S(CN_{\ge -n}(A)[-n])\ar[r]_I\ar[u]^1&S(CP_{\ge -n}(A)[-n]).\ar[u]_{\iota}}
\]
The lemma is immediate from this.
\end{proof}

If we identify both rows in the lemma above and use diagram \eqref{fundechar}, we obtain a commutative diagram with
exact rows
\begin{equation}\label{htpyfundechar}
\xymatrix{ KH_{n+1}A\ar[r]\ar[d]_{ch_{n+1}}& K^{\nil}_nA\ar[r]\ar[d]_{\nu_n}&K_n A\ar[r]\ar[d]_{c_n}& KH_nA
\ar[d]_{ch_n}\ar[r]&\tau K^{\nil}_{n-1}A\ar[d]_{\nu_{n-1}}\\
               HP_{n+1}A\ar[r]_S&HC_{n-1}A\ar[r]_B&HN_nA\ar[r]_I&HP_nA\ar[r]_S&HC_{n-2}A.}
\end{equation}

We call $\nu_n:K^{\nil}_nA\to HC_{n-1}A$ the {\it secondary (Chern) character} and $ch_n:KH_nA\to HP_nA$ the {\it homotopy
Chern character}.

\subsection{Comparison with the secondary character of M. Karoubi and C. Weibel}\label{subsec:compakw}
The construction presented above is a nonconnective variant of Weibel's construction \cite{wenil};
the latter predates the invention of $KH$ and uses Karoubi-Villamayor $K$-theory instead. Recall the
Karoubi-Villamayor $K$-theory space is $KV(A)=K(\DA A)$. Weibel considers the map of spaces
$c:K(A)\to SCN_{\ge 1}A$, applies it to $\DA A$, and uses the equivalences
\[
SCN_{\ge 1}(\DA A)\weq SCP_{\ge 1}(\DA A)\overset\sim\leftarrow SCP_{\ge 1}A
\]
to obtain a homotopy commutative diagram
\[
\xymatrix{K^{\nil'}A\ar[d]_{\nu'}\ar[r]&KA\ar[d]_c\ar[r]&KV(A)\ar[d]_{ch'}\\
          S(CC(A)[-1])\ar[r]_(0.56)B&SCN_{\ge 1}A\ar[r]&SCP_{\ge 1}A.}
\]
There is a map $j$ from $KV(A)$ to the connected component of ${}_0\KH(A)$, which is not a weak equivalence in general.
However if $A$ is $K_0$-regular, $j$ is an equivalence (\cite{wh}), whence $KV_nA\cong KH_nA$ for $n\ge 1$, and
$K^{\nil}_nA=K^{\nil'}_nA$ for $n\ge 1$ and is zero for $n\le 0$. Moreover it is clear that in this case
$ch'_n=ch_n$ and $\nu'_n=\nu_n$ for $n\ge 1$.

Weibel's $\nu'$ is related to a relative character introduced by Karoubi \cite{karast} in the context of Fr\'echet
algebras which is discussed in Section \ref{sec:locass}. According to \cite[3.4]{karmult}, the algebraic version of
Karoubi's relative character coincides with Weibel's. By \cite[6.17]{karast}, $K^{\nil'}_*(A)=\pi_*(GL(\DA A)/GL(A))^+$. Similarly,
$HC_{*-1}(A)=\pi_*(SN_{\ge 1}(\DA A)/SN_{\ge 1}(A))$ and $\nu'$ is the map induced by $ch'$.

\subsection{Extension to nonunital algebras}
Let $A\in\ass$, a not necessarily unital algebra. Write $\tilde{A}_\Q=A\oplus \Q$ for its unitalization. If $X$ is a
functor from $\ass_1$ to spectra, define a functor $Y$ from $\ass$ to spectra by
\[
YA:=\hofi(X\tilde{A}_\Q\to X\Q).
\]
If $A$ happens to be unital then $\tilde{A}_\Q\cong A\times \Q$, whence if $X$ preserves products up to homotopy,
we have $X(\tilde{A}_\Q)\weq XA\times X\Q$ and $XA\weq YA$. Thus any such functor $X$ can be extended
up to homotopy to all of $\ass$; we shall abuse notation and write $XA$ for $YA$ ($A\in\ass$). With this convention we
have $XA\weq X(\tilde{A}_\Q:A)$ for all $A\in\ass$. In particular this applies
to all functors appearing in diagram \eqref{fundechar}, and so the latter makes sense and is homotopy commutative
for all $A\in\ass$. Hence the primary, homotopy and secondary characters are defined for all $A\in\ass$. Alternatively
one can perform the same construction using $\tilde{A}_\Z$ instead of $\tilde{A}_\Q$. The following lemma shows that
this change leads to equivalent spectra.

\begin{lem}
Let $A$ be a $\Q$-algebra, $X(A)$ any of the spectra in diagram \eqref{fundechar}. Then the natural map
$X(\tilde{A}_\Z:A)\to X(\tilde{A}_\Q:A)$ is a weak equivalence.
\end{lem}
\begin{proof}
Note that
\[
C_p(\tilde{A}_\Z:A)=A^{\otimes p}\oplus A^{\otimes p+1}=C_p(\tilde{A}_\Q:A)
\]
whence $C(\tilde{A}_\Z:A)=C(\tilde{A}_\Q:A)$. Since each of $CP_n$, $CC_n$ and $CN_n$ is a product of copies of some terms of $C$, it
follows that $CP(\tilde{A}_\Z:A)=CP(\tilde{A}_\Q:A)$, $CC(\tilde{A}_\Z:A)=CC(\tilde{A}_\Q:A)$ and
$CN(\tilde{A}_\Z:A)=CN(\tilde{A}_\Q:A)$. By \ref{hntauk} and the proof of \ref{htodotautodo}, each of the level spaces
of each of the spectra of the bottom row of \eqref{fundechar} is weak equivalent to the space the Dold-Kan correspondence assigns to a
truncation of one of the cyclic complexes. Thus the assertion of the lemma is true for each of the spectra in the bottom
row of \eqref{fundechar}. To prove it is also true for $\K$, we must show that
\begin{equation}\label{flit}
\K(\tilde{A}_\Z,\tilde{A}_\Q:A)=\hofi(\K(\tilde{A}_\Z:A)\to\K(\tilde{A}_\Q:A))
\end{equation}
is weakly contractible. Because $A$ is a $\Q$-algebra, it is excisive for $K$-theory with finite coefficients (\cite[1.9]{qs},
\cite{chu}).
Hence it suffices to show that the rational homotopy groups
\[
K_*^\Q(\tilde{A}_\Z,\tilde{A}_\Q:A)=\pi_*\K(\tilde{A}_\Z,\tilde{A}_\Q:A)\otimes\Q
\]
vanish. By \cite{kabi},
\begin{align*}
K_*^\Q(\tilde{A}_\Z,\tilde{A}_\Q:A)&=HC_{*-1}(\tilde{A}_\Z,\tilde{A}_\Q:A)\otimes_\Z\Q\\
                                   &=HC_{*-1}(\tilde{A}_\Z\otimes\Q,\tilde{A}_\Q:A)\\
                                   &=HC_{*-1}(\tilde{A}_\Q,\tilde{A}_\Q:A)=0.
\end{align*}
Thus the lemma holds for $\K$. A similar argument shows that
\[
\K(\DA_n(\tilde{A}_\Z),\DA_n(\tilde{A}_\Q):\DA_nA)\weq 0
\]
for each $n$, whence $\K(\DA(\tilde{A}_\Z),\DA(\tilde{A}_\Q):\DA A)\weq 0$. Thus the lemma
is also true for $\KH$. The assertion for the remaining spectra follow from the fact that each
of the rows and columns in \eqref{fundechar} is a homotopy fibration.
\end{proof}

\subsection{An application: excision for $KH$}
Weibel has shown that $KH$ satisfies excision in the category of all rings \cite{wh}; a different proof is given in \cite[8.1.2]{wt}.
As an application of the material discussed in this section, we give yet another proof for $\Q$-algebras.

\begin{propo}\label{khex}
Let
\begin{equation}\label{ses}
0\to A'\to A\to A''\to 0
\end{equation}
be an exact sequence in $\ass$. Then the sequence
\begin{equation}\label{khfib}
\KH A'\to \KH A\to \KH A''
\end{equation}
is a homotopy fibration.
\end{propo}
\begin{proof}
By \eqref{fundechar}, we have a functorial homotopy fibration
\[
\K^{\inf}(\DA(?))\to \KH(?)\to \tau\K(\DA(?)).
\]
By Lemma \ref{htodotautodo} and the excision property of $HP$, $\tau\K(\DA(?))$ maps \eqref{ses} to a homotopy fibration.
Because $\DA$ is exact and $\K^{\inf}$ excisive, the same is true of $\K^{\inf}(\DA(?))$. It follows that \eqref{khfib}
is a homotopy fibration.
\end{proof}

\section{$K^{\inf}$-regularity}\label{sec:kinfreg}
\subsection{Vorst's method}
Let $F$ be a functor from algebras to abelian groups. Recall a ring $R$ is called $F$-regular if $F(R)\to F(R[t_1,\dots,t_m])$
is an isomorphism for every $m\ge 1$. In this section we investigate $K^{\inf}_n$-regularity. We say that $R$ is
{\it $K^{\inf}$-regular} if it is $K^{\inf}_n$-regular for all $n$. The first basic result we shall establish is an analogue of Vorst's theorem \cite{vorst}. The method of proof follows Vorst's.
We need some notations and preliminaries. As is costumary, if $X$ is a functor from rings to spectra we write $NX(A)$
for the fiber of the map $XA[t]\to X A$ induced by $t\mapsto 0$.
Let $A$ be a unital algebra, $A[t]$ the polynomial ring, $f\in A$ a central element,
$A_f:=A[f^{-1}]$ the localized algebra, and
$\eta_f:A[t]\to A[t]$ the $A$-algebra homomorphism determined  by $t\mapsto ft$. Consider the direct system $\{A_n\}_{n\ge 1}$
with $A_n=A[t]$ and structure maps $\eta_f:A_n\to A_{n+1}$ for all $n$. We have $\coli_n A_n=A+tA_f[t]\subset A_f[t]$.
Put $K^{\inf}A_{[f]}=\coli_n K^{\inf}A_n$. We have a natural map
\begin{equation}\label{voma}
NK^{\inf}A_{[f]}\to NK^{\inf}(A_f).
\end{equation}
\begin{lem}\label{lem:coli}
The functors $N_mK_n^{\inf}$ preserve filtering colimits.
\end{lem}
\begin{proof}
In view of \eqref{fundechar}, we have $N_mK_n^{\inf}=N_mK_n^{\inf,nil}$, and a long exact sequence involving $N_mK_*^{\inf,nil}$ as
well as $N_mK^{nil}_*=N_mK_*$ and $N_mHC_*$. The lemma follows from this and the fact that $K_*$ and $HC_*$ preserve filtering colimits.
\end{proof}
\begin{lem}\label{localiza}
Let $A$ be a unital algebra, $f\in A$ a central element which is not a zero divisor. Then the map \eqref{voma} is
a weak equivalence.
\end{lem}
\begin{proof}
Vorst proved the analogous statement for $K$-theory in \cite[Lemma 1.4]{vorst}, using that $K$-theory preserves filtering
colimits. The same argument, using \ref{lem:coli}, proves the lemma.
\end{proof}

\begin{propo}\label{propo:vokinf}
Let $A\in\ass$, $n\in\Z$. Assume $A$ is $K^{\inf}_n$-regular. Then $A$ is $K^{\inf}_{n-1}$-regular.
\end{propo}
\begin{proof}
Write $\tilde{A}:=\tilde{A}_\Q$. Consider the commutative diagram
\[
\xymatrix{\tilde{A}\ar[d]\ar[r]&\tilde{A}[t^{-1}]\ar[d]\\
          \tilde{A}\ar[r]&\tilde{A}[t,t^{-1}].}
\]
Because both $\K$ and (in view of \ref{hntauk}) $\tau\K$ send this diagram to a homotopy cartesian diagram, so
do $\K^{\inf}$ and $N\K^{\inf}$. In particular we have a commutative diagram with exact rows
\[
\xymatrix{0\to NK_n^{\inf}\tilde{A}\ar[d]\ar[r]&NK_n^{\inf}\tilde{A}[t]\oplus NK^{\inf}_n\tilde{A}[t^{-1}]\ar[r]\ar[d]&
NK^{\inf}_n\tilde{A}[t,t^{-1}]\ar[r]\ar[d]&NK_{n-1}^{\inf}\tilde{A}\to 0\ar[d]\\
0\to NK_n^{\inf}{\Q}\ar[r]&NK_n^{\inf}{\Q}[t]\oplus NK^{\inf}_n{\Q}[t^{-1}]\ar[r]&
NK^{\inf}_n{\Q}[t,t^{-1}]\ar[r]&NK_{n-1}^{\inf}\Q\to 0.}
\]
We point out for later use that because all vertical arrows are surjective, the sequence of their kernels is exact.
The hypothesis of the proposition says that the first vertical map is an isomorphism. By Lemma \ref{localiza} this
implies that also the arrow next to the rightmost one is an isomorphism. In particular its kernel vanishes, whence
the same must be true of that of the rightmost arrow.
\end{proof}

\begin{propo}\label{propo:nuiso}
Let $A\in\ass$. If $A$ is $K_n^{\inf}$ regular, then the secondary character $\nu_m:K^{\nil}_mA\to HC_{m-1}A$ is an isomorphism
for $m\le n$.
\end{propo}
\begin{proof}
Because $\K^{\inf}\DA A$ is defined as the total spectrum of a simplicial spectrum, there is a spectral sequence
\[
E^2_{p,q}=\pi_p([r]\mapsto K^{\inf}_q(\DA_r A))\Rightarrow K^{\inf}_{p+q}(\DA A).
\]
By \ref{propo:vokinf}, the hypothesis implies that, for $q\le m$, $[r]\mapsto K_q(\DA_r A)$ is a constant simplicial abelian group.
Hence the map $K^{\inf}_mA\to K^{\inf}_m(\DA A)$ is an isomorphism for $m\le n$ and a surjection for $m=n+1$. It follows
that $K_m^{\inf,\nil}A=0$ for $m\le n$, whence $\nu_m$ is an isomorphism for $m\le n$.
\end{proof}
\begin{coro}\label{coro:rel_nuiso}
Let $A\to B$ be a surjective homomorphism with kernel $I$. Assume that both $A$ and $B$ are $K^{\inf}_n$-regular. Then the
secondary character induces an isomorphism $\nu_m:K^{\nil}_m(A:I)\to HC_{m-1}(A:I)$ for $m<n$.
\end{coro}
\begin{proof}
Immediate from \eqref{propo:nuiso} and the $5$-lemma.
\end{proof}
\subsection{$K^{\inf}$-regularity vs. $K$-regularity}\label{subsec:regkink}
\begin{lem}\label{regu} Let $A$ be a $\Q$-algebra.
\begin{enumerate}
\item[(i)]If $n\le -1$ then $A$ is $K_n^{\inf}$-regular if and only if it is $K_n$-regular.
\item[(ii)]\label{k0infreg} If $A$ is $K_0^{\inf}$-regular then it is $K_0$-regular.
\item[(iii)] Any nonzero commutative unital noetherian regular $A\in\ass_1$ is $K$-regular and not $K_0^{\inf}$-regular.
\end{enumerate}
\end{lem}
\begin{proof}
By \eqref{basicfib} we have an exact sequence
\[
NHN_{n+1}A\to NK^{\inf}_nA\to NK_nA\to NHN_nA.
\]
If $n\le 0$ then $NHN_{n}A=0$; this proves (i) and (ii). On the other hand we have
\[
NHN_{1}A=NHC_0A=\bigoplus_{n\ge 1}\frac{A}{[A,A]}t^n.
\]
If $A$ is unital noetherian regular then it is $K$-regular, whence $NK_0^{\inf}A=NHN_1(A)$,
which is nonzero if $A\ne 0$ is commutative. This proves (iii).
\end{proof}
\begin{coro}\label{khkv} If $A$ is $K_0^{\inf}$-regular, then $KH_n(A)=K_n(A)$ if $n\le 0$ and
$KH_n(A)=KV_n(A)$ if $n\ge 1$.
\end{coro}
\subsection{Examples}\label{subsec:exaregkinf}
\begin{propo}\label{propo:notcom}
Let $A$ be a unital algebra, $[A,A]\subset A$ the commutator subspace, $\langle[A,A]\rangle$ the two-sided ideal it generates.
Assume $A$ is $K^{\inf}_n$-regular for some $n\ge 1$. Then $A=\langle[A,A]\rangle$.
\end{propo}
\begin{proof}
Assume the inclusion $\langle[A,A]\rangle\subset A$ is strict. Then $B=A/\langle[A,A]\rangle$ is a nonzero commutative unital algebra, and therefore
it surjects onto a field $F\supset \Q$. Because fields are $K$-regular, $K^{\nil}_*F=0$, whence the composite of
$\nu_*:K^{\nil}_*A\to HC_{*-1}A$ with $HC_{*-1}A\to HC_{*-1}F$ is zero. If $A$ is $K^{\inf}_n$-regular and $n\ge 1$,
this implies that $HC_{0}A\to HC_{0}F$ is zero, by \ref{propo:nuiso}.
In particular the natural inclusion $\Q=HC_0\Q\to HC_0F=F$ must be zero, which is absurd.
\end{proof}
\begin{exa}\label{exa:nilpo}
Assume there is an $n$ such that $A^n=0$. Then $(\DA_mA)^n=0$ for all $m$, and because $K^{\inf}$ is nilinvariant,
 $K^{\inf}_*(A)=K^{\inf}_*(\DA A)=0$. In particular $A$ is $K^{\inf}$-regular.
\end{exa}
\begin{lem}\label{promv}
Let
\[
\xymatrix{A\ar[d]\ar[r]&B\ar[d]\\
          C\ar[r]&D}
\]
be a pullback diagram in $\ass$ such that the vertical arrows are onto. If any three of $A$, $B$, $C$ and $D$ are
$K^{\inf}$-regular, then so is the fourth.
\end{lem}
\begin{proof}
Tensoring the diagram above with $\Q[t]$ we still get a pullback diagram with surjective vertical maps.
Because $\K^{\inf}$ is excisive, it maps both diagrams to a homotopy cartesian diagrams of spectra. Hence $NK^{\inf}$
also sends it to a homotopy cartesian diagram. The lemma is immediate from the Meyer-Vietoris sequence
\[
NK_{n+1}^{\inf}D\to NK_n^{\inf}A\to NK_n^{\inf}B\oplus NK_n^{\inf}C\to NK_n^{\inf}D\to NK_{n-1}^{\inf}A
\]
\end{proof}

\begin{coro}\label{coronilpo} Assume $A\triqui B$ is an ideal in an algebra such that $B^n\subset A$ for some $n\ge 1$. Then $A$ is
$K^{\inf}$-regular if and only if $B$ is.
\end{coro}
\begin{proof} Immediate from \ref{exa:nilpo} and \ref{promv}.
\end{proof}

Another source of examples is provided by Wagoner's {\it infinite sum algebras} \cite[2.4]{wa}. Recall from \cite{wa} that
a {\it sum algebra} is a unital algebra $A$ together with elements $\alpha_i,\beta_i$, $i=0,1$ such that the following
identities hold
\begin{gather*}
\alpha_0\beta_0=\alpha_1\beta_1=1\\
\beta_0\alpha_0+\beta_1\alpha_1=1
\end{gather*}
If $A$ is a sum algebra, then
\begin{gather*}
\oplus:A\times A\to A,\\
(a,b)\mapsto a\oplus b=\beta_0a\alpha_0+\beta_1b\alpha_1
\end{gather*}
is a homomorphism in $\ass_1$. An {\it infinite sum algebra} is a sum algebra together with a unit preserving
algebra homomorphism $\infty:A\to A$, $a\mapsto a^\infty$ such that
\begin{equation}\label{ainfi}
a\oplus a^\infty=a^\infty\qquad (a\in A).
\end{equation}
If $A$ is a sum algebra and $A\hookrightarrow C$ a monomorphism in $\ass_1$, then $C$ is again a sum algebra, with the same
choice of $\alpha_i,\beta_i$. In particular if $B$ is any unital algebra, then $A=A\otimes 1\subset C=A\otimes B$, whence
$A\otimes B$ is a sum algebra; if moreover $A$ is an infinite sum algebra then $A\otimes B$ is one too, with $\infty$-map
$(a\otimes b)^\infty=a^\infty\otimes b$. Applying this to $B=\Q[t],\Q[t,t^{-1}], M_n\Q$, we get that $A[t]$, $A[t,t^{-1}]$
and $M_nA$ are (infinite) sum algebras if $A$ is.

\begin{lem}\label{lem:infisum} Let $A\in\ass_1$ be an infinite sum algebra. Then for every $B\in \ass$ the algebra $A\otimes B$ is
$K^{\inf}$-regular and $K^{\inf}_*(A\otimes B)=0$.
\end{lem}
\begin{proof}
Put $\tilde{B}=\tilde{B}_\Q$. By the discussion above, $A\otimes \tilde{B}$ is an infinite sum algebra. Hence if the
lemma is true for $B=\Q$, it is true for all $B\in\ass$, by \ref{promv} and excision. The proof of \cite[6.4.1]{wt} shows that
if $A$ is a sum algebra and $H_*$ is any $M_2$-stable, excisive homology theory, then $H_*(A)=0$. Applying this to $K^{\inf}_*$
we obtain the lemma.
\end{proof}

\section{Locally convex algebras} \label{sec:locass}
\subsection{Preliminaries}

In this paper a {\it locally convex vectorspace} is a $\C$-vectorspace with a complete locally
convex topology. We write $\hotimes$ for the completed projective tensor  (over $\C$) of locally convex
vectorspaces. The category $\LCV$ of locally convex vectorspaces and continuous linear maps, equipped with $\hotimes$, is a
symmetric monoidal category. A {\it locally convex algebra} is a monoid object in this monoidal category (not necessarily with
unit). We write
$\locass$ for the category of locally convex algebras, and $\locass_1$ for the subcategory of unital algebras and unit
preserving maps. Substituting $\hotimes$ for $\otimes$ in the definition of the cyclic complexes, we obtain
complexes $C^{\top}$, $CC^{\top}$, etc, whose homology we denote by $HH^{\top}$, $HC^{\top}$, etc. A short exact sequence
\[
0\to A'\to A\to A''\to 0
\]
in $\locass$ is called {\it admissible} if the forgetful functor $\locass\to \LCV$ maps it to a split exact sequence.
A functor $X$ from $\locass$ to spectra is called {\it excisive} if it maps admissible sequences to homotopy fibrations.
It is proved in \cite{val} (see also \cite{cv}) that the $HP^{\top}$ satisfies excision, or, in the notation we have
just introduced, that the spectrum $\HP^{\top}$ which defines it is excisive. Put

\[
\A^n:=\{(x_0,\dots,x_n)\in\R^{n+1}:\sum x_i =1 \}\supset \Delta^n:=\{x\in \A^n:x_i\ge 0 \ (0\le i\le n)\}.
\]

If $V\in\LCV$ we write
\[
\DD_n V:=C^\infty(\Delta^n) \hotimes V.
\]
Here, $C^\infty(\Delta^n)$ denotes the locally convex vectorspace of all those functions on $\Delta^n$ which are restrictions
of $C^\infty$-functions $\A^n\to \C$.
The cosimplicial structure on $[n]\mapsto\Delta^n$ induces
a simplicial one on $\DD V$. In particular, for $A\in\locass$, $\DD A$ is a simplicial locally convex algebra.
Let $d_0,d_1:\DD_1 A\to A$ be the face maps. Two homomorphisms $f_0,f_1:B\to A\in \locass$ are {\it diffeotopic} if
there exists $H:B\to \DD_1A\in \locass$ such that the following diagram commutes
\[
\xymatrix{&\DD_1 A\ar[d]^{(d_0,d_1)}\\
          B\ar[ur]^H\ar[r]_{(f_0,f_1)}&A.}
\]

A functor $X$ from $\locass$ to spectra is {\it diffeotopy invariant} if it sends diffeotopic maps to homotopic maps.
Equivalently $X$ is diffeotopy invariant if for every $A\in \locass$ the image $X(s_0)$ of the degeneracy $A\to \DD_1A$
is a weak equivalence. The maps
\begin{gather*}
h_i:\DD_n\DD_1 A\to \DD_{n+1},\\ h_i(f)(t_0,\dots,t_n,x)=f(t_0,\dots,t_i+t_{i+1},\dots, t_n,(t_{i+1}+\dots t_n)x)
\end{gather*}
$0\le i\le n$ form a simplicial homotopy between the identity and the map $\DD_n(s_0d_0)$. Hence if $X$ is any functor
from $\locass$ to spectra, then the total fibrant spectrum $X\DD A$ of the simpilicial spectrum $[n]\mapsto X\DD_nA$ is
diffeotopy invariant. The spectral sequence of a simplicial spectrum gives
\[
E^2_{p,q}=\pi_q([r]\mapsto \pi_pX(\DD_rA))\Rightarrow \pi_{p+q}X(\DD A).
\]
If $X$ happens to be diffeotopy invariant already, then $[r]\mapsto\pi_pX(\DD_rA)$ is constant, whence the natural
map $XA\to X(\DD A)$ is a weak equivalence in this case. Thus $X$ is diffeotopy invariant if and only if $X\to X\DD$
is a weak equivalence. For later use we define
\[
X^{{\rm rel}}A=\hofi(XA\to X^{\dif}A).
\]

\begin{lem}\label{contrapi}
Let $V\in\LCV$; then the Moore complex $\M\DD V$ has a continuous linear contracting homotopy.
\end{lem}
\begin{proof}
Formula \eqref{hcontra} defines the required homotopy.
\end{proof}

\begin{coro}\label{bufo}
 If $A\in\locass$ then  $HC_*(\DD A)=HH_*(\DD A)=0$ and $HN_*(\DD A)\cong HP_*(\DD A)$.
\end{coro}
\begin{proof}
Follows from lemmas \ref{pinulhcnul} and \ref{contrapi}.
\end{proof}

\begin{defi} Let $A\in\locass$. The {\it diffeotopy $K$-theory spectrum} is $\KD(A):=\K(\DD A)$.
\end{defi}

\begin{propo}\label{propo:kdexc}
$\KD$ is diffeotopy invariant and excisive.
\end{propo}
\begin{proof}
As explained above, $A\mapsto X\DD A$ is diffeotopy invariant for any funtorial spectrum $X$, in particular
for $X=\K$. The proof of excision for $\KD$ is analogous to that for $\KH$ (\ref{khex}); simply substitute \ref{bufo}
for \ref{bafa}, and recall that $\hotimes$ preserves admissible exact sequences.
\end{proof}

\subsection{Primary and secondary characters}\label{pschdif}

\begin{lem}\label{kunetop} Let $V$, $W$ be simplicial objects in $\LCV$, $V\hotimes W:[n]\mapsto V_n\hotimes W_n$ their completed projective
tensor product. Assume $\M V$ has a continuous linear contracting homotopy. Then $\pi_*(V\hotimes W)=0$.
\end{lem}
\begin{proof}
By the Eilenberg-Zilber theorem, we have a spectral sequence
\[
E^1_{p,q}=H_p(\M(V\hotimes W_q))\Rightarrow \pi_{p+q}(V\hotimes W)
\]
On the other hand it is clear that $H_*((\M V)\hotimes W_q)=0$  under the hypothesis of the lemma.
Hence it suffices to show that $\M(V\hotimes W_q)=\M V\hotimes W_q$. By the Dold-Kan theorem, if $U$ is
any simplicial abelian group, then $U=S\M U$, where $S$ is a functor sending each nonnegatively graded chain complex $C$ to
a simplicial abelian group $SC$ with
\[
SC_n=\bigoplus_{p=0}^n C_p\otimes_\Z \Lambda ^p\Z^n
\]
Specifically, for $\{e_1,\dots,e_n\}$ is the canonical basis of $\Z^n$ and $i_1<\dots<i_p$, the isomorphism $S\M U\iso U$ sends
\[
\M_pU\otimes\Lambda^p\Z^n\owns x\otimes v_{i_1}\land\dots\land v_{i_p}\mapsto s_{i_p}\circ\dots\circ s_{i_1}(x)\in U_n
\]
If $U$ is a simplicial locally convex space, then the image of each of the terms in the decomposition above is a closed subspace, because
$ s_{i_p}\circ\dots\circ s_{i_1}$ is a section with continuous right inverse $d_{i_1}\circ\dots\circ d_{i_p}$. It follows
that the decomposition is preserved upon tensoring with a constant simplicial locally convex space. In particular
$\M(V\hotimes W_q)=\M V\hotimes W_q$ as wanted.
\end{proof}

\begin{lem}\label{isotopo} The maps $HN^{\top}_*\DD A\to HP^{\top}_*\DD A\leftarrow HP^{\top}_*A$ are isomorphisms.
\end{lem}
\begin{proof}
It follows from \ref{kunetop} that $HH_*(\DD A)=HC_*(\DD A)=0$, whence the first map is an isomorphism.
That the second map is an isomorphism is a consequence of the diffeotopy invariance of $HP^{\top}$ \cite{val}.
\end{proof}

Next we define a variant of the characters of Section \ref{pschalg} for $A\in\locass$. First of all, from Section
\ref{pschalg} we obtain a commutative diagram with exact rows
\[
\xymatrix{KD_{n+1}A\ar[d]_{c_{n+1}\DD}\ar[r]&K^{\rel}_nA\ar[d]_{c^{\rel}_n}\ar[r]&K_nA\ar[d]_{c_n}\ar[r]&KD_nA\ar[d]_{c_n\DD}\\
           HN_{n+1}\DD A\ar[r]_S&HN^{\rel}_{n}A\ar[r]_B&HN_n A\ar[r]_I&HN_n\DD A}
\]
Composing with the natural map $HN_*\to HN^{\top}_*$ and using the isomorphisms of Lemma \ref{isotopo} we obtain
a commutative diagram
\[
\xymatrix{KD_{n+1}A\ar[d]_{ch^{\dif}_{n+1}}\ar[r]&K^{\rel}_nA\ar[d]_{ch^{\rel}_n}\ar[r]&K_nA\ar[d]_{\hat{c}_n}\ar[r]&KD_nA\ar[d]_{ch^{\dif}_n}\\           HP^{\top}_{n+1}A\ar[r]_S&HC^{\top}_{n-1}A\ar[r]_B&HN^{\top}_n A\ar[r]_I&HP_n^{\top}A}
\]
Here $ch^{\dif}$ stands for {\it diffeotopy Chern character} and $ch^{\rel}$ for {\it relative Chern character}.

\subsection{Comparison with the algebraic characters of Section \ref{pschalg}}

The map $\hat{c}_*:K_*A\to HN_*^{\top}A$ has been defined by the commutativity of the diagram
\[
\xymatrix{K_*(A)\ar[r]^{\hat{c}_*}\ar[d]_{c_*}&HN_*^{\top}A\\
           HN_*A\ar[ur]&}
\]
To compare the remaining algebraic and topological characters we need some preliminaries.
 Let $X$ be a functor from $\locass$ to spectra and $A\in\locass$. Consider the natural maps $\iota: A\to \DA A$, $j:A\to\DD A$.
We have a homotopy commutative diagram
\[
\xymatrix{XA\ar[d]_{X(j)}\ar[rr]^{X\iota}&&X\DA A\ar[d]^{X(\DA(j))}\\
          X\DD A\ar[rr]_{X(\iota\DD)}&&X\DA\DD A}
\]

\begin{lem} The map $X(\iota\DD)$ is a weak equivalence.
\end{lem}
\begin{proof}
The multiplication map of any countably dimensional $\C$-algebra $\Lambda$ is continuous in the fine topology, whence every
such algebra can be considered as an locally convex algebra. Moreover we have $\Lambda\hotimes V=\Lambda\otimes_\C V$ for any
$V\in\LCV$. In particular this applies to $\Lambda=\DA_n\C$, so $\DA_nA\in\locass$ for every $n\ge 0$ and every $A\in\locass$.
Moreover
\[
\DD_p(\DA_q A)=\DD_p\C\hotimes\DA_q\C\hotimes A=\DA_q\otimes\DD_pA=\DA_q(\DD_pA).
\]
Thus it suffices to show that $A\mapsto X\DD A$ is invariant under polynomial homotopy. Because it is already diffeotopy
invariant, is enough to show that the map $s_0d_0:\DA_1A\to \DA_1A$ is diffeotopic to the identity. One checks that
$H:\DA_1A\to \DD_1\DA_1A$, $H(f)(s,t)=f(st)$ is a homotopy $s_0d_0\to 1$.
\end{proof}

Applying the lemma to $X=\K$, $\tau\K$ and $\HN$ and using the isomorphisms of Lemmas \ref{bafa} and \ref{isotopo},
we obtain commutative diagrams
\begin{equation}\label{compa_diags}
\xymatrix{KH_*A\ar[d]_{ch_*}\ar[r]&KD_*A\ar[d]^{ch^{\dif}_*}\\
          HP_*A\ar[r]& HP^{\top}_*A}
\qquad\xymatrix{K^{\nil}_*A\ar[d]_{\nu_*}\ar[r]&K^{\rel}_*A\ar[d]^{ch_*^{\rel}}\\
                HC_{*-1}A\ar[r]&HC_{*-1}^{\top}A.}
\end{equation}
In particular we have a comparison map
\begin{equation}\label{compa_map}
KH_*(A)\to KD_*(A).
\end{equation}
\subsection{Comparison with Karoubi's relative character}\label{subsec:compaktop}
As in the polynomial case, one can consider, for any $A\in\locass$, the diffeotopic version
$KV^{\dif}A:=K(\DD A)$ of Karoubi-Villamayor's K-theory. For $A$ a Fr\'echet algebra, $KV^{\dif}A$ is Karoubi's topological $K$-theory
(\cite{karast}, \cite{karmult}). He puts $K^{\rel'}A=\hofi(KA\to KV^{\dif}A)$ and defines a relative character $ch^{\rel'}_n:K^{\rel'}_nA\to HC_{n-1}^{\top}A$ by means
of a rather explicit integral formula (see \cite{karast}). In \cite{karmult} he shows that $ch^{\rel'}_*$ fits into a commutative diagram
\[
\xymatrix{KV^{\dif}_{n+1}\ar[d]_{ch^{\dif'}_{n+1}}\ar[r]&K^{\nil'}_nA\ar[d]_{ch^{\rel'}_n}\ar[r]&K_nA\ar[d]_{\hat c_n}\ar[r]&
KV^{\dif}_n(A)\ar[d]_{ch^{\dif'}_n}
\\ HP^{\top}_{n+1}A\ar[r]_S&HC^{\top}_{n-1}A\ar[r]_B&HN^{\top}_nA\ar[r]_I&HP_n^{\top}A}
\]
constructed in a similar way as in the algebraic case considered by Weibel. One checks that if the canonical map
\[
K_pA\to K_p\DD_qA
\]
is an isomorphism for all $p\le 0$ and all $q\ge 0$, then for $n\ge 1$, $KV^{\dif}_nA=KD_nA$, $K^{\rel'}_n A=K^{\rel}_nA$,
$ch^{\rel'}_n=ch^{\rel}_n$, and $ch^{\dif'}_n=ch^{\dif}_n$; for $n\le 0$ we have $ch^{\rel}_n=0$.

\section{Operator ideals and stability}\label{sec:opid_stab}

\subsection{Operator ideals}\label{subsec:opideals}
An {\it operator ideal} is a functor from the category $\hilb$ of infinite-dimensional, separable, complex Hilbert spaces
and isometries to the category of $\C$-algebras, together with a natural ideal embedding $\J(H)\triqui \B(H)$ into
the algebra of bounded operators. We shall abuse notation and write $\J$ and $\B$ for $\J(H)$ and $\B(H)$. We implicitly
assume that $\J\ne 0$ and $\J \ne \B(H)$.

\begin{rem} Any ideal of the algebra of bounded operators of any particular Hilbert space $H_0$ can be made into
an operator ideal in the sense of the definition above; see \cite[3.3]{hus}. We will use this fact without mentioning.
\end{rem}

\begin{nota}\label{nota:opid}
We say that an operator ideal $\J$ is {\it locally convex} if it carries a complete locally convex topology which makes it into
a locally convex algebra and is such that the inclusion $\J\to\B$ is continuous. The operator ideal $\J$ is a {\it Fr\'echet ideal}
if it carries a complete metrizable locally convex topology such that the inclusion
$\J\to \B$ is continuous. This implies, as we shall see presently, that the multiplication map $\B \times \J \times \B \to \J$ is
jointly continuous; in particular $\J$ is locally convex. Indeed, consider the map of left multiplication $\B \times \J \to \J$.
If we fix $a \in \B$, then the map $b\mapsto ab$ is continuous from $\J$ to $\B$ and hence to $\J$, by the closed graph theorem which applies
to linear maps between Frech\'et spaces, see \cite[Thm. $2.15$]{Rud}.
Similarily, for $b \in \J$ the map $a \mapsto ab$ is continuous from $\B$ to $\B$ and hence to $\J$, again, using the closed graph
theorem. This shows that the map $\B \times \J \to \J$ is separately continuous and hence jointly continuous, using that all algebras
involved are Fr\'echet algebras, see \cite[Thm. $2.17$]{Rud}.
A similar remark applies to the multiplication map $$\B \times \J \times \B \to \J.$$

A Fr\'echet ideal $\J$ is a {\it Banach ideal} if it is a Banach space. Note that our definition of Banach ideal has built-in the condition
that its topology can be defined by a norm $||\ \ ||_\J$ satisfying
\[
||abc||_\J\le||a||\cdot ||b||_\J\cdot||c||
\]
for all $a,c\in\B(H)$ and $b\in \J$.

We use the following notations for operator ideals. We write $\cL_p$ for the $p$-Schatten ideal ($p>0$), $\cL_\infty$ for the ideal of all compact operators,
and $\F$ for that of finite rank operators.

Write $\circtimes$ for the completed tensor product of Hilbert spaces. Let
$\I,\J$ be operator ideals.
Consider the natural map $\boxtimes:\B(H)\otimes\B(H)\to \B(H\circtimes H)$. We say that $\I$ {\it multiplies} $\J$ if
\[
\I(H)\boxtimes \J(H)\subset \J(H\circtimes H).
\]
An operator ideal $\J$ is {\it multiplicative} if it multiplies itself.
\end{nota}
\sn
\begin{propo}\label{propo:subharm}
Let $H\in\hilb$, $\B:=\B(H)$, and $\J \subset \B$ a locally convex ideal
in the sense of \ref{nota:opid} above. Assume that the multiplication map
\[
\B\times \J \times \B \to \J
\]
is jointly continuous. Then
\begin{enumerate}
\item[(i)] $\cL_1\subset\J$.
\item[(ii)]  $\cL_1$ multiplies $\J$.
\end{enumerate}
\end{propo}
\begin{proof}
Since $\J$ is an operator ideal, all finite rank operators belong to $\J$. In particular, the rank-one partial isometries
$\phi_{\xi,\eta} = \xi \langle\eta,.\rangle$ belong to $\J$. Let $\rho$ be a continuous seminorm on $\J$.
Since $\B\times \J \times \B\to \J$ is jointly continuous, there is a continuous
seminorm $\sigma$ such that
\[
\rho(abc) \leq C \|a\| \sigma(b) \|c\|,\quad \forall a,c \in \B, b \in \J.
\]
Let $\xi,\eta \in H$ be unit-vectors and let $\alpha \in H$ be a fixed unit-vector. We have
\begin{equation}\label{rule}
\phi_{\xi,\eta}=\phi_{\xi,\alpha}\circ\phi_{\alpha,\alpha}\circ\phi_{\alpha,\xi}
\end{equation}
Thus
\[
\rho(\phi_{\xi,\eta}) \leq C \|\phi_{\xi,\alpha}\| \sigma(\phi_{\alpha,\alpha}) \| \phi_{\alpha,\eta}\| \leq C_{\rho},
\]
for some constant $C_{\rho}$ and all unit-vectors $\xi,\eta \in H$.
By the spectral theorem, every operator $a \in \cL_1$ has a presentation:
\begin{equation}\label{spectral_rep}
a = \sum_{i=1}^{\infty} \lambda_i \phi_{\xi_i,\eta_i}
\end{equation}
with $(\lambda_i)_{i \in \N}$ summable. From the boundedness of $\rho(\phi_{\xi,\eta})$ by $C_{\rho}$,
it follows that the sequence of partial sums $a_n = \sum_{i=1}^n \lambda_i \phi_{\xi_i,\eta_i}$
is a Cauchy sequence with respect to $\rho$. Since $\J$ is complete and $\rho$
arbitrary, this implies that $a_n$ converges in $\J$. Since $\J$ as well as $\cL_1$ sit continuously in $\B$,
the limit has to be $a \in \cL_1$. We conclude that $\cL_1 \subset \J$. This proves i). To prove ii), consider
\[
T(\J) = \{ a \in \B(H)\colon a \boxtimes b \in \J(H \circtimes H),\, \forall b \in \J(H)\}.
\]
It is clear that $T(\J)$ is an ideal, and that all finite rank operators belong to $T(\J)$.
We want to show that $\cL_1 \subset T(\J)$. Let $a\in \cL_1$ be as in \eqref{spectral_rep} and $b\in \J$. Using \eqref{rule},
and proceeding as above, one shows that the sequence of partial sums
\[
a_n \boxtimes b = \sum_{i=1}^n  \lambda_i \phi_{\xi_i,\eta_i} \boxtimes b
\]
converges in $\J(H\circtimes H)$ to the operator $a\boxtimes b$. This finishes the proof.
\end{proof}
\begin{exa}\label{exa:frech_subhar}
Any Fr\'echet operator ideal $\J$ satisfies the hypothesis of \ref{propo:subharm}. In particular $\cL_1\subset \J$ and $\cL_1$ multiplies $\J$.
Moreover, as also $\cL_1$ is Fr\'echet, the map $\boxtimes:\cL_1(H)\otimes \J(H)\to \J(H\circtimes H)$ is jointly continuous, whence it extends
to a continuous homomorphism
\[
\boxtimes:\cL_1(H)\hotimes\J(H)\to \J(H\circtimes H).
\]
\end{exa}

\subsection{Periodicity}\label{subsec:perio}
Write $\com$ for the algebra of smooth compact operators (\cite[1.4]{cd}), $\T_0$ for the pointed Toeplitz algebra \cite[6.4]{cd}, and
$B[0,1]$ the sub-algebra of $C^{\infty}(\Delta^1) \hotimes B$ formed by functions
all of whose derivatives vanish at the endpoints $0,1$. Similarly, put
$B[0,1) = \ker(ev_1\colon B[0,1] \to B)$ and $B(0,1)= \ker(ev_0 \oplus ev_1\colon B[0,1] \to B \oplus B)$.
For any $B\in\locass$ we have admissible exact sequences
\begin{equation}\label{admiseq}
\xymatrix{0\ar[r]&B\hotimes\com\ar[r]&B\hotimes\T_0\ar[r]&B(0,1)\ar[r]&0\\
          0\ar[r]&B(0,1)\ar[r]&B[0,1)\ar[r]&B\ar[r]&0}
\end{equation}
Let now $X$ be an
excisive, diffeotopy invariant functor from $\locass$ to the derived category of spectra. Using excision and the fact
that $B[0,1)$ is diffeotopically contractible, we obtain a natural weak equivalence
\[
\Omega XB\iso XB(0,1).
\]
In the next proposition and elsewhere we consider, for $A\in\ass$, the canonical inclusion
\begin{equation}\label{e11}
e_{11}:A\to M_2A=A\otimes_\C M_2\C,\quad a\mapsto e_{11}(a):=a\otimes e_{11}
\end{equation}
A functor $X$ from $\locass$ to the derived category of spectra is called {\it $M_2$-stable} if it maps \eqref{e11} to a weak equivalence.
\begin{propo}\label{propo:periodicity}
Let $X$ be a functor from $\locass$ to the derived category of spectra and $\J$ a Fr\'echet operator ideal.
Assume $X$ is excisive, diffeotopy invariant, and $M_2$-stable. There is a natural isomorphism $X(A\hotimes\J)\weq\Omega^2 X(A\hotimes \J)$.
\end{propo}
\begin{proof}
From Proposition \ref{propo:subharm} we have $\com \subset\cL_1\subset\J$.
Note that the inclusion $\com \subset \J$ is continuous. It is shown in the proof of \cite[Prop. 6.1.2]{ct} that if $\J$ is multiplicative the
inclusion $M_2\C\to \com$ induces a diffeotopy equivalence
$j: \J\otimes_\C M_2\C\to \J \hotimes\com$. Actually the same argument applies using only that $\J$ is multiplied by $\cL_1$, which is always
the case, by \ref{propo:subharm} ii). Thus $X$ is $\com$-stable, since it is $M_2$-stable and diffeotopy invariant.
It is proved in \cite[Satz 6.4]{cd} that if $F:\locass\to\ab$ is a functor which is diffeotopy invariant, split exact and $\com$-stable,
then $F$ satisfies $F(A\hotimes\T_0)=0$. Actually in {\it loc.\,cit.} this result
is stated only for functors defined on the subcategory of $m$-algebras, but the proof holds for $\locass$ as well.
 In particular this applies to each of the functors
$\pi_n(X(? \hotimes\J)):\locass\to \ab$, ($n\in\Z$). Thus $X(A\otimes\T_0)=0$.
The proposition is now immediate from the sequences \eqref{admiseq}.
\end{proof}

\begin{coro}\label{perikd} Let $A$ be a locally convex algebra and let $\J$ be a Fr\'echet operator ideal as above.
There is a natural isomorphism $KD_*(A\hotimes \J) \cong KD_{*+2}(A\hotimes \J)$.
\end{coro}
\begin{proof}
Since we already know from \ref{propo:kdexc} that $\KD$ is diffeotopy invariant and excisive, we just have to check that
it is $M_2$-stable. This means that for all $A\in\locass$, $\KD$ applied to \eqref{e11} is an equivalence. Because $\KD$ is
excisive, it suffices to show this for $A\in\locass_1$.
But this follows from the fact that $\K$ maps \eqref{e11} to
a weak equivalence for every $A\in\ass_1$.
\end{proof}

\subsection{Bivariant $K$-theory}\label{subsec:kk}
Bivariant topological $K$-theory for locally convex algebras (defined by Cuntz in \cite{cw}) associates to any pair $(A,B)$ of such
algebras a $\Z$-graded abelian group $kk^{\lc}_*(A,B)$. The graded groups $kk^{\top}_*$ are periodic of period $2$, satisfy excision in both variables, are diffeotopy invariant
and moreover they are $\com$-stable, in the sense that the natural inclusion $A\to A\hotimes \com$ induces isomorphisms
$kk^{\lc}_*(A\hotimes\com,?)\weq kk^{\lc}_*(A,?)$,
and similarly in the other variable. In addition $kk_*^{\lc}$ is equipped with an associative composition product
\begin{equation}\label{compo_top}
kk^{\lc}_p(B,C)\otimes kk^{\lc}_q(A,B)\to kk^{\lc}_{p+q}(A,C).
\end{equation}
In particular one can consider a category $kk^\lc$ whose objects are those of $\locass$ and where the homomorphisms are given by
$\hom_{kk^\lc}(A,B)=kk^\lc_0(A,B)$.
Cuntz shows that $kk$ is universal in the sense that every functor $E:\locass\to\ab$ which is diffeotopy invariant, half exact and $\com$-stable
extends uniquely to $kk$. In addition to composition, there is defined an external or cup product
\begin{equation}\label{cup_lc}
kk^{\lc}_p(A,B)\otimes kk^{\lc}_q(C,D)\to kk^{\lc}_{p+q}(A\hotimes B,C\hotimes D)
\end{equation}
satisfying the obvious associativity rules, and compatible with the composition product above.

The groups $kk_*^\lc(A,B)$ are difficult to compute in general. For example, not much is known of $kk_0(\C,\C)$, apart
from the fact that it is nonzero (see \cite{cw}). On the other hand, it is shown in \cite{ct} that if $\J$ is a harmonic operator ideal
(see \ref{defi:harm} below for a definition), then
\begin{equation}\label{kk0=k0}
kk_0(\C,A\hotimes \J)=K_0(A\hotimes \J).
\end{equation}
An algebraic version of Cuntz' theory, called $kk$, was introduced in \cite{wt}. Algebraic $kk$ is defined for arbitrary rings, is equipped
with a composition product, and satisfies a similar universal property as $kk^{\top}$;  simply substitute polynomial homotopy invariant for diffeotopy invariant and $M_\infty$-stable for $\com$-stable.
It is proved in \cite{wt} that if $A$ is any ring, then there is a natural isomorphism
\[
kk_*(\Z,A)=KH_*(A).
\]
Note the ring $\Z$ appears instead of $\C$; this is because we are considering arbitrary rings, i.e. $\Z$-algebras.
Alternatively, the universal property mentioned above can be formulated for algebras over any fixed ring $\ell$,
and excision considered only with respect to exact sequences which are split in some fixed underlying category,
such as those of sets, abelian groups, $\ell$-modules or
$\ell\otimes_\Z\ell^{op}$-modules. Each of these choices leads to a different variant of $kk$. However if $\ell$ is commutative, and $kk^{\ell}_*$ is
the variant relative to the category of $\ell$-algebras which are central as $\ell$-bimodules (e.g. all $\C$-algebras considered in this paper are
central in this sense) and to excision with respect to any of the choices of exact sequences mentioned above, then we have
\begin{equation}\label{kk=kh}
kk_*^{\ell}(\ell,A)=KH_*(A).
\end{equation}
For general $\ell$ there is no analogue of the cup product \eqref{cup_lc}. However if $\ell$ is a field,
we do have an associative cup product (\cite[6.6.5]{cw})
\begin{equation}\label{cup_ell}
kk^{\ell}_p(A,B)\otimes kk^{\ell}_q(C,D)\to kk^{\ell}_{p+q}(A\otimes_\ell B,C\otimes_\ell D).
\end{equation}
In particular this gives a cup product for $KH$, using \eqref{kk=kh}.

\section{Main results}\label{sec:main}
\subsection{Diffeotopy invariance}

\begin{defi}\label{defi:harm}
Let $\{e_i:i\in \N\}$ be the canonical Hilbert basis of $l^2(\N)$. The {\it harmonic operator} $\omega\in\B(l^2(\N))$
is the following diagonal operator:
\[
\omega(e_n)=\frac{e_n}{n}.
\]
A {\it harmonic ideal} is a multiplicative Banach operator ideal $\J$ such that $\omega \in \J(l^2(\N))$.
\end{defi}
\sn
The following result, due to Cuntz and the second author, will be used extensively in what follows. The formulation we present here is not
given in \cite{ct}, but follows immediately from \cite[Thms. 4.2.1, 5.1.2]{ct}.
\sn
\begin{thm}[Diffeotopy Invariance Theorem] \label{hit} Let $E:\locass\to\ab$ be a split exact, $M_2$-stable functor, and $\J$ a harmonic ideal.
Then the functor $A\to E(A\hotimes \J)$ is diffeotopy invariant.
\end{thm}
\begin{coro}\label{corohit}
Let $X$ be a functor from $\locass$ to the derived category of spectra. If $X$ is split-exact and $M_2$-stable, and $\J$ is a harmonic Banach ideal,
then $A\mapsto X(A\hotimes \J)$ is diffeotopy invariant.
\end{coro}

By definition, a harmonic ideal is necessarily Banach. We shall need a variant of \ref{hit}
which is valid for all Fr\'echet ideals $\J$. First we need some notation. Let $\alpha:A\to B$
be a homomorphism of locally convex algebras. We say that $\alpha$ is an {\it isomorphism up to square zero}
if there exists a continous linear map $\beta:B\hotimes B\to A$ such that the compositions $\beta\circ(\alpha\hotimes\alpha)$
and $\alpha\circ\beta$ are the multiplication maps of $A$ and $B$. Note that if $\alpha$ is an isomorphism up to square zero, then
its image is a ideal of $B$, and both its kernel and its cokernel are square-zero algebras.

\begin{defi}\label{defi:nilinv} Let $E:\locass\to\ab$ be a functor. We call $E$ {\rm continously nilinvariant} if it sends isomorphisms
up to square zero into isomorphisms.
\end{defi}

\begin{exa}\label{exa:cont_nilinv} It is proved in \cite[Thm. 2.3.3]{ct} that any isomorphism up to square zero $\alpha$ defines an invertible element in
$kk^{\top}$. It follows from this and the universal property of $kk^{\top}$, that any excisive, diffeotopy invariant, $\com$-stable
functor $E:\locass\to\ab$, sends $\alpha$ to an isomorphism. Thus such functors are continuously nilinvariant. For another
example, assume $E$ is the restriction to $\locass$ of a functor $\ass\to\ab$ which is half exact and nilinvariant.
Then $E$ is continuously nilinvariant.
\end{exa}

\begin{thm}\label{thm:frech_hit}
Let $E:\locass\to\ab$ be a split exact, $M_2$-stable and continuously nilinvariant functor, and $\J$ a Fr\'echet ideal.
Then the functor $A\to E(A\hotimes \J)$ is diffeotopy invariant.
\end{thm}
\begin{proof}
By \ref{hit}, the theorem is true for $\J=\cL_2$. By continous nilinvariance, the inclusion $\cL_1\to\cL_2$ induces a natural
isomorphism $E(?\hotimes\cL_1)\to E(?\hotimes\cL_2)$. In particular the theorem holds for $\J=\cL_1$. Now let $\J$ be an arbitrary
Fr\'echet ideal, and $H\in\hilb$. Choose an isometry $H\circtimes H\cong H$. The composite
\[
\xymatrix{\phi:\J(H)\ar[rr]^(0.4){1_{\J}\otimes e_{11}}&&\J(H)\hotimes\cL_1(H)\ar[r]^{\boxtimes}&\J(H\circtimes H)\ar[r]^(.6){\cong}&\J(H)}
\]
is a morphism of the form $a\mapsto UaU^*$ for some partial isometry $U$. Because $E$ is $M_2$-stable and split exact, it maps
$1_A\otimes \phi:A\hotimes\J\to A\hotimes \J$ to an isomorphism, by \cite[Lemma 3.2.3 (c)]{ct}. Thus $F(?):=E(?\hotimes\J)$ is naturally
a retract of $G(?):=E(?\hotimes\J\hotimes\cL_1)$, which is diffeotopy invariant, as we have just shown. But then $F$ must
map the inclusion $A\to A[0,1]$ to an isomorphism, because $G$ does, since the retract of an isomorphism is an isomorphism.
This shows that $F$ is diffeotopy invariant, concluding the proof.
\end{proof}
\begin{coro}\label{coro:frech_hit}
Let $X$ be a functor from $\locass$ to the derived category of spectra. If $X$ is excisive, $M_2$-stable and continously nilinvariant,
and $\J$ is a Fr\'echet operator ideal, then $A\mapsto X(A\hotimes \J)$ is diffeotopy invariant.
\end{coro}

\subsection{Fr\'echet ideals and $K^{\inf}$-regularity}\label{subsec:main}

\begin{thm}\label{thm:main}
Let $L\in\locass$, $\J$ a Fr\'echet operator ideal
and $A$ a $\C$-algebra.
\begin{enumerate}
\item[(i)] The functors $\K^{\inf}(A\otimes_\C(?\hotimes\J))$ and $\KH(A\otimes_\C(?\hotimes\J))$ are diffeotopy invariant.
\item[(ii)] The algebra $A\otimes_\C( L\hotimes\J)$ is $K^{\inf}$-regular.
\item[(iii)] The comparison maps
{\rm \eqref{compa_map}} and
{\rm \eqref{kk0=k0}} induce natural isomorphisms
\[
KD_n(L\hotimes \J)=kk^{\top}_n(\C,L\hotimes \J)=KH_n(L\hotimes \J) \quad \forall n\in\Z.
\]
Moreover we have
\[
KH_n(L\hotimes \J)=KV_n(L\hotimes \J)=KV^{\dif}_n(L\hotimes \J) , \quad n\ge 1.
\]
\end{enumerate}
\end{thm}
\begin{proof} Part i) is immediate from \ref{coro:frech_hit} and the nilinvariance, excision, and $M_2$-stability properties of
$\K^{\inf}$ and $\KH$.
Part ii) follows from i), since polynomially homotopic maps are diffeotopic. From \ref{coro:frech_hit} and the nilinvariance,
excision, and $M_2$-stability properties
of $\KH$, we get that $\KH(?\hotimes\J)$ is diffeotopy invariant. It follows that the comparison map \eqref{compa_map} is an isomorphism.
Next note that by i) and \ref{regu}, $L\hotimes\J$ is $K_0$-regular. In particular, by \eqref{kk0=k0}, we have natural isomorphisms
\begin{equation}\label{kk0=khj}
KH_0(L\hotimes\J)=K_0(L\hotimes\J)=kk_0^{\top}(\C,L\hotimes\J)
\end{equation}
As both $KH_*(?\hotimes\J)$ and $kk^{top}_*(\C,?\hotimes\J)$ are excisive, $\com$-stable and diffeotopy invariant,
\eqref{kk0=khj} implies that they are isomorphic. This proves the first assertion of iii). The second follows from the
$K_0$-regularity of $L\hotimes\J$ and the fact, proved above, that  $K_{-n}(?\hotimes\J)=KH_{-n}(?\hotimes\J)$ is diffeotopy invariant.
\end{proof}
\begin{rem}
Note that the preceding result implies that $K_0(?\hotimes\J)$ is diffeotopy invariant for every Fr\'echet ideal $\J$. This improves
the results concerning diffeotopy invariance of $K_0$ that
were obtained in \cite{ct}. There, Theorem \ref{hit} was applied to show that $K_0(? \hotimes \J)$ satisfies diffeotopy invariance for
$\J$ harmonic, and in particular for $\J=\cL_p$, $p>1$.
The result above implies that also $K_0(? \hotimes \cL_1)$ satisfies diffeotopy invariance. This refinement
is possible through the detailed analysis of $K^{\inf}$-regularity and was beyond the scope of the methods
used in \cite{ct}.
\end{rem}
\subsection{Cyclic homology and the comparison map $K\to K^{\top}$}\label{subsec:cyclic}
Theorem \ref{thm:main} implies that the various possible definitions of covariant topological $K$-theory coincide after stabilizing by
Fr\'echet ideals. We unify notation and put
\[
K_*^{\top}(L\hotimes \J):=KH_*(L\hotimes \J)=KD_*(L\hotimes \J)=kk^{\top}_*(\C,L\hotimes \J).
\]

\begin{thm} \label{thm:exact} Let $L \in \locass$ and let $\J$ be a Fr\'echet operator ideal. For each $n \in \Z$,
there is
a natural $6$-term exact sequence of abelian groups as follows:
\[
\xymatrix{ K^{\top}_{1}(L\hotimes \J)\ar[r]&HC_{2n-1}(L\hotimes \J)\ar[r]&K_{2n}(L\hotimes \J)\ar[d] \\
K_{2n-1}(L \hotimes \J) \ar[u]& HC_{2n-2}(L \hotimes \J) \ar[l] & K_0^{\top}(L\hotimes \J). \ar[l]}
\]
\end{thm}
\begin{proof} Immediate from \ref{propo:nuiso}, \eqref{htpyfundechar}, and Theorem \ref{thm:main}.
\end{proof}

\begin{rem}\label{rem:chrel-nu}
It follows from Theorem \ref{thm:main} and Subsection \ref{subsec:compaktop} that if $\J$ is a Fr\'echet operator ideal, then
$K_*^{\rel}(\J)=K_*^{\nil}(\J)$ and Karoubi's relative character agrees with the character
$ch^{rel}_*:K^{\rel}_*(\J)\to HC_{*-1}^{\top}(\J)$ defined in
Subsection \ref{pschdif}. In particular, by \eqref{compa_diags} it factors through the algebraic character
$\nu_*:K_*^{\nil}(\J)\to HC_{*-1}(\J)$.
\end{rem}

\subsection{$K_0$ and $K_{-1}$ for operator ideals}\label{subsec:k0k-1}

The result of this subsection, due to M. Wodzicki, is stated without proof in \cite{wod}. We give a short and complete proof.
\begin{thm} \label{thm:compk0k1}
Let $H\in\hilb$ and let $\J \subset \B(H) $ be an operator ideal. We have natural isomorphisms
\[ K_0(\J) = \Z \qquad \mbox{and} \qquad K_{-1}(\J) = 0.\]
\end{thm}
\begin{proof}
It was proved in \cite{karcomp} that $K_0(\cL_{\infty}) = \Z$ and
$K_{-1}(\cL_{\infty})=0$. Note this also follows from Theorem \ref{thm:main} and Karoubi's density theorem \cite[Thm 1.4]{karcomp}.
Next we consider the $\C$-unitalization of various subalgebras of $\B$. In this proof, the superscript $?^+$ will
denote $\C$-unitalization; thus $A^+=A\oplus \C$. If $A\subset \cL_\infty$, then $A^+$ identifies with the subalgebra $A\oplus \C$ of $\B$.
Consider the extension
\begin{equation}\label{compext}
0 \to \J \to \cL_{\infty}^+ \to \cL_{\infty}^+/\J \to 0
\end{equation}
The $K$-theory sequence associated to \eqref{compext}, together with the augmentation map
$(?)^+\to \C$, give a commutative diagram with exact row
\begin{equation}\label{k0k-1}
\xymatrix{K_1(\cL_{\infty}^+)\ar[r]^{\pi_1}&K_1(\cL_{\infty}^+/\J)\ar[r]^{\alpha}&K_0(\J)\ar[r]&K_0(\cL_{\infty}^+)\ar[dr]\ar[r]^{\pi_0}&
K_0(\cL_{\infty}^+/\J)\ar[d]\ar[r]&K_{-1}(\J)\to 0.\\
          && & &\Z &}
\end{equation}
Next note that, because for all $k \in \cL_{\infty}$ and $\lambda \in \C^\times$,
the operator $\lambda1 + k$ is invertible in $\cL_\infty^+/\F$, the ring $\cL_\infty^+/\J$ is local.
Thus the vertical map in \eqref{k0k-1} which involves $K_0$ is an isomorphism (\cite[1.3.11]{rosen}). On the other hand,
$K_0(\cL_\infty^+)\to \Z=K_0(\C)$ is surjective, as it
is induced by the augmentation $\cL_\infty^+\to \C$. It follows that $\pi_0$ is onto, and thus $K_{-1}(\J)=0$ as claimed.
It remains to show that $\pi_1$ is surjective.
Since $\cL_\infty^+/\J$ is local, $K_1(\cL_\infty^+/\J)$ identifies with $(\cL_\infty/\J+\C^\times)_{ab}$ (\cite[2.2.6]{rosen}).
Now any element of $(\cL_\infty/\J+\C^\times)_{ab}$
lifts in $\cL_\infty^+$ to a sum $T=k+\lambda$, with $k$ a compact operator and $\lambda\in\C^\times$. Such an operator $T$ is clearly
Fredholm and of Fredholm index zero. Thus we can find a finite rank partial isometry $f$ which maps $\ker(T)$ isomorphically
onto $\ran(T)^\perp$, so that $T+f$ is invertible in  $\cL_\infty^+$. It follows that $\alpha$ is zero, whence
$K_0(\J)\cong \ker(K_0(\cL_\infty^+)\to \Z)=K_0(\cL_\infty)=\Z$.
\end{proof}

\subsection{Sub-harmonic ideals and $K^{\inf}$-regularity}\label{subsec:main_alg}
In this subsection we prove a variant of Theorem \ref{thm:main} (see \ref{thm:main_alg}) which is valid in a more algebraic context.
We need some notation. Let $\J$ be an operator ideal. The {\it root completion} of $\J$ is the ideal
\[
\J_{\infty} = \bigcup_{n\geq 1} \J^{\nicefrac1n}.
\]
Thus $\J_{\infty}$ is an operator ideal and $(\J_{\infty})^2 = \J_{\infty}$. Note that,
contrary to what our notation might suggest:
\[
(\cL_1)_{\infty} \neq \cL_{\infty},
\]
despite the fact that $(\cL_1)^{\nicefrac1p} = \cL_{p}$ for all $p>0$.

\begin{defi}\label{defi:subhar}
Let $\J$ be an operator ideal. We say that $\J$ is {\it sub-harmonic} if there exists $p>0$ such that
$\cL_p$ mutliplies $\J_\infty$.
\end{defi}

\begin{exa} Every Fr\'echet ideal is sub-harmonic, by \ref{propo:subharm} ii). The same is true of every Schatten ideal $\cL_p$
($p>0$), since they are all multiplicative.
\end{exa}

\begin{thm}\label{thm:main_alg}
Let $A$ be a $\C$-algebra and $\J$ a sub-harmonic ideal. Under these hypothesis, we have:
\begin{enumerate}
\item[(i)] The algebra $A\otimes_\C\J$ is $K^{\inf}$-regular.
\item[(ii)] The graded abelian group $KH_*(A \otimes_\C\J)$ is $2$-periodic; we have
\begin{equation}\label{formu_perio}
KH_n(A \otimes_\C\J)=\left\{\begin{matrix} K_0(A\otimes_\C\J)& \text{ if }n\text{ is even.}\\
                                           K_{-1}(A\otimes_\C\J)&\text{ if }n\text{ is odd.}\end{matrix}\right.
\end{equation}
\end{enumerate}
\end{thm}
\begin{proof}
Applying part i) of \ref{thm:main} with $\J=\cL_1$ and $L=\C$, we obtain part i) of the current theorem for $\cL_1$. By the
nilinvariance and excision properties of $K^{\inf}$, we get that if $p<q<\infty$ then the inclusion $\cL_p\subset\cL_q$ induces
a natural isomorphism $K^{\inf}_*(A\otimes\cL_p)=K^{\inf}_*(A\otimes\cL_q)$. It follows that i) is true for $\J=\cL_p$ ($p>0$).
Now let $\J$ be an operator ideal and $p>0$ such that $\cL_p$ multiplies $\J$. Essentially same argument of the proof of \ref{thm:frech_hit}
(one just needs to substitute $\cL_p$ and $\otimes_\C$ for $\cL_1$ and $\hotimes$) shows that
$K^{\inf}_*(?\otimes\J)$ is naturally a retract of $K^{\inf}_*(?\otimes\J\otimes\cL_p)$. But we have just seen that $A\otimes\J\otimes\cL_p$
is $K^{\inf}$-regular. Thus i) holds whenever $\cL_p$ multiplies $\J$. In particular if $\cL_p$ multiplies $\J_\infty$ then
$A\otimes\J_\infty$ is $K^{\inf}$-regular. But by \ref{lem:coli} and \ref{coronilpo}, this implies that also $A\otimes\J$ is $K^{\inf}$-regular,
proving i). To prove ii) if suffices to consider the case $\J=\J_\infty$, since all functors involved are nilinvariant and preserve
filtering colimits, and since excision holds for both $KH$ and nonpositive $K$-theory. Thus we may assume that $\cL_p$ multiplies
$\J$. The multiplication maps $\cL_p(H)\otimes\cL_p(H)\to\cL_p(H\circtimes H)\cong \cL_p(H)$, together with the cup product \eqref{cup_ell},
make $KH_*(\cL_p)$ into a ring. Further, the map $\cL_p(H)\otimes\J(H)\to\J(H\circtimes H)\cong \J(H)$ defines a $KH_*(\cL_p)$-module structure.
By \ref{thm:main}, the nilinvariance and excision properties of $KH$, and \ref{thm:compk0k1}, we have
\[
KH_*(\cL_p)=KH_*(\cL_1)=K^{\top}_*(\cL_1)=\Z[u,u^{-1}].
\]
Here $u$ is homogeneous of degree $2$. This proves that $KH_*(A\otimes_\C\J)$ is $2$-periodic. Formula \eqref{formu_perio} follows
from this and part i).
\end{proof}

\section{Some applications}\label{sec:appli}

\subsection{Wodzicki's exact sequence}
In this subsection we prove a generalized version of a $6$-term exact sequence involving relative $K$-theory and cyclic homology due to M.
Wodzicki (\cite[Theorem 5]{wod}). His result is that the sequence \eqref{relseq} in Theorem \ref{thm:absolute_wodzicki} below is exact for $H$-unital $A$.
We extend this to all algebras, and give a version involving absolute, rather than relative $K$-theory and cyclic homology.
\sn
Wodzicki's result appeared without proof in \cite[Theorem 5]{wod}. He has lectured in several places giving full details of his
proof; in \cite{hus}, D. Husem\"oller took the effort to write down the details for the case $A=\C$. Our proof of
\ref{thm:absolute_wodzicki}
is different as it uses Theorem \ref{thm:main_alg}, which in turn is based on Theorems \ref{thm:main} and \ref{thm:frech_hit}, and ultimately on the homotopy invariance theorem \ref{hit} and on the excision theorem for $K^{\inf}$ from \cite{kabi}, none of which were available at the time when \cite{wod} was written.
\sn
\begin{thm} \label{thm:absolute_wodzicki}
Let $\J$ be a sub-harmonic ideal and let $A$ be a $\C$-algebra. For each $n \in \Z$,
there exist two natural $6$-term exact sequences of abelian groups as follows:
\begin{gather}
\label{absoseq}\xymatrix{ K_{-1} (A\otimes_\C \J)\ar[r]&HC_{2n-1}(A\otimes_\C \J)\ar[r]&K_{2n}(A\otimes_\C \J)\ar[d] \\
K_{2n-1}(A \otimes_\C \J) \ar[u] & HC_{2n-2}(A \otimes_\C \J) \ar[l] & K_{0}(A\otimes_\C \J). \ar[l]}\\
\label{relseq}\xymatrix{ K_{-1} (A\otimes_\C \J)\ar[r]&HC_{2n-1}(A\otimes_\C\B:A\otimes_\C \J)\ar[r]&K_{2n}(A\otimes_\C\B:A\otimes_\C \J)\ar[d] \\
K_{2n-1}(A\otimes_\C\B:A\otimes_\C \J) \ar[u] & HC_{2n-2}(A\otimes_\C\B:A \otimes_\C \J) \ar[l] & K_{0}(A\otimes_\C \J). \ar[l]}
\end{gather}
\end{thm}
\begin{proof} The assertions are immediate from Theorem \ref{thm:main_alg}, Proposition \ref{propo:nuiso}, and Corollary
\ref{coro:rel_nuiso}.
\end{proof}
\subsection{More on $K$-theory of operator ideals}
The following propositions subsume the results of \ref{thm:absolute_wodzicki} and \ref{subsec:main_alg} for the special case $A=\C$.
\begin{propo}\label{propo:kopi} Let $\I$ be a sub-harmonic ideal, $n\in\Z$, and $p\ge 1$. Then
\begin{enumerate}
\item[(i)] $KH_{2n}(\I)=\Z$ and $KH_{2n+1}(\I)=0$.
\item[(ii)] There is a natural exact sequence
\[
0\to HC_{2n-1}(\I)\to K_{2n}(\I)\to \Z\overset{\alpha_n}\to HC_{2n-2}(\I)\to K_{2n-1}(\I)\to 0.
\]
Here, under the identification $KH_{2n}(\I)=\Z$ of {\rm \eqref{formu_perio} and \ref{thm:compk0k1}}, the map $\alpha_n$
identifies with the composite $S\circ ch_{2n}$ of \eqref{htpyfundechar}.
\item[(iii)] If $\I\subset \cL_p$, then $\alpha_n$ is injective for all $n\ge (p+1)/2$.
\end{enumerate}
\end{propo}
\begin{proof}
Assertion i) follows from Theorem \ref{thm:compk0k1} in combination with Theorem \ref{thm:main_alg}.
Now ii) is immediate from i) and Theorem \ref{thm:absolute_wodzicki}.
To prove iii) it suffices to consider the case $\I=\cL_p$. By Remark \ref{rem:chrel-nu}, we can use \cite[4.13]{karop}; the
assertion is immediate from this.
\end{proof}

\begin{exa} As a particular case of \ref{propo:kopi} iii), we obtain that if $\I\subset \cL_1$ is sub-harmonic then
\[
K_1(\I)=(\I/[\I,\I]) / (\Z S(ch_2(\xi_2)))
\]
where $\xi_2$ is a generator of $KH_2(\I)=\Z$. For example this applies when $\I=\cL_\epsilon$, $0<\epsilon\le 1$.
\end{exa}

\begin{propo}\label{propo:k1opi}
Let $\I$ be a Fr\'echet operator ideal and assume that $[\I,\I]$ contains some positive operator.
\begin{enumerate}
\item[(i)] There are natural isomorphisms
\[ K_0(\I) = \Z \quad  \mbox{and} \quad K_1(\I) = \I/[\I,\I]. \]
\item[(ii)] There is an exact sequence
\[
\xymatrix{ 0 \ar[r] & HC_1(\I) \ar[r] & K_2(\I) \ar[r] & \Z \ar[r] & 0.  }
\]
\end{enumerate}
\end{propo}
\begin{proof} If $[\I,\I]$ contains some positive operator, then by Theorem $5.11$ of \cite{dykwodz}, it contains the whole principal ideal
which is generated by this operator. In particular, all operators of finite rank $\F$ are in $[\I,\I]$.
The inclusion $\C\cong \C e_{11}\subset\F\subset \I$ induces an isomorphism $\Z=K_2^{\top}(\C)\cong K_2^{\top}(\I)$. Thus the image
of $S\circ ch^{\dif}_2:KH_2(\I)\cong K_2^{\top}(\I)\to HC_0^{\top}(\I)=\I/[\I,\I]=HC_0(\I)$ is contained in the image of
$\C e_{11}$ in $\I/[\I,\I]$. By our observation $ \F \subset [\I,\I]$ and hence the map $\alpha_1$ has to be zero.
\end{proof}

\begin{rem}\label{rem:k1harmonic}
Let $\J$ be a harmonic Banach ideal. In what follows, we use the computation
\[
K_1(\J)=\J/\J^2
\]
of \cite[7.1]{ct}. We point out that there is a mistake in {\it loc.\,cit.};
the proof given there needs that not only $\J$ but also $\J^2$ be harmonic. This is necessary to guarantee that the element $h$ considered
in \cite[pp. 370, line 11]{ct} actually is in $\J^2$. Furthermore we shall see in Remark \ref{rem:pearcy-topping} below that the
conclusion of \cite[7.1]{ct} does not hold for $\J=\cL_2$.
\end{rem}

The following proposition follows from \cite[Thm.\ 5.11]{dykwodz}, specifically from the implication ii)(c)$\Rightarrow$ii)(a)
in {\it loc.\ cit.} We want to point out that in \cite{dykwodz} a complete structure theory for operator ideals is developed, and much more general results can be obtained. However, using our techniques, we can give a different proof of a particular case of interest. 

\begin{propo}\label{propo:k1harm}
Let $\J$ be a Banach operator ideal such that $\J^2$ is harmonic. Then
\[
[\J,\J]=\J^2.
\]
\end{propo}
\begin{proof} Consider the ring homomorphism $\J\to \J/\J^2$. By \ref{propo:kopi} and \ref{exa:nilpo},
we have a commutative diagram with exact rows
\[
\xymatrix{\Z\ar[r]^{\alpha_1}\ar[d]&HC_0(\J)\ar[d]\ar[r]&K_1(\J)\ar[d]\ar[r]&0\\
          0\ar[r]& HC_0(\J/\J^2)\ar[r]_{\cong}&K_1(\J/\J^2)\ar[r]&0}
\]
Consider the determinant map $\det:K_1(\J/\J^2)\to (1+\J/\J^2)^*\cong \J/\J^2$. From the description of the relative character given in
\ref{subsec:compakw}
it is straightforward to check that $\det$ is the inverse of the isomorphism of the second row. It was proved in \cite[7.1]{ct}
that the composite of $\det$ with the second vertical map is an isomorphism. It follows that ${\rm Im}(\alpha_1) \cong \J^2/[\J,\J]$. Since
the right hand side is a vectorspace and the left hand side a quotient of $\Z$, we conclude that $\alpha_1=0$. This implies the claim.
\end{proof}
\begin{rem}\label{rem:pearcy-topping}
Proposition \ref{propo:k1harm} applies in particular if $\J=\cL_p$ and $p>2$, showing that $[\cL_p,\cL_p] = \cL_{p/2}$, which is a classical result of Pearcy and Topping,
see \cite{PT}. The methods which were used in \cite{ct} are very much inspired by those of \cite{PT}. It is therefore not surprising
that one of the main results of \cite{PT} can be derived as a corollary to the results obtained in \cite{ct}. On the other hand,
since $\omega^2\in \cL_1=\cL_2^2$ but the arithmetic mean of its eigenvalues is not summable, we
have $\omega^2\notin [\cL_2,\cL_2]$, by \cite[5.6]{dykwodz}. Thus the inclusion $[\cL_2,\cL_2]\subset\cL_1$ is strict. Note that the
argument of the proof of \ref{propo:k1harm} above shows that if $\J$ is any sub-harmonic ideal such that the
determinant map $K_1(\J)\to \J/\J^2$ is an isomorphism, then $\J^2=[\J,\J]$. It follows that
the map $det:K_1(\cL_2)\to\cL_2/\cL_1$ is surjective but not injective. Thus
\[
\ker(det:K_1(\cL_2)\to\cL_2/\cL_1)\ne 0
\]

Note also that for $p>2$, we obtain
\begin{equation}\label{k1lp}
K_1(\cL_p)=\cL_p/\cL_{p/2}\ne 0\ \ (p>2).
\end{equation}
This observation was made by the first author in 2004, using Pearcy-Topping's result and the diagram considered in the
proof of \ref{propo:k1harm} above. As this contradicted the well-established calculation $K_1(\cL_p)=0$ of
\cite[4.1]{karop}, finding a direct proof of \eqref{k1lp} seemed in order. This was done by the second author (see \cite[\S7]{ct}).
\end{rem}

\subsection{Application to the representation algebras of A. Connes and M. Karoubi}\label{subsec:ck}

Throughout this section, $\I$ will be a fixed operator ideal.
Let $H\in\hilb$ be a Hilbert space; put $\B=\B(H)$. Consider the following two subalgebras
of $\B(H\oplus H)$.
\begin{gather*}
\M^0(\I)=\left\{\left[\begin{matrix}u&0\\0&v\end{matrix}\right]: u-v\in \I\right\}\cong \B\times_{\B/\I}\B\\
\M^1(\I)=\left[\begin{matrix}\B&\I\\ \I&\B\end{matrix}\right]
\end{gather*}
These algebras were considered in \cite{ck} in the case when $\I=\cL_p$ and $1\le p\le\infty$; for $1\le p<\infty$ and
$j\equiv p-1{\rm\ \ mod} (2)$, they are the universal
algebras for $p$-summable Fredholm modules.

\begin{propo}\label{propo:mji} Let $\I$ be an operator ideal and $j=0,1$.
\begin{enumerate}
\item[(i)] If $H_*$ is an excisive, $M_2$-invariant homology theory, then $H_*(\M^j(\I))=H_{*+j}(\I)$.
\item[(ii)] If $\I$ is $K^{\inf}$-regular, then so is $\M^j(\I)$.
\item[(iii)] If $\I$ is sub-harmonic, and $n\in\Z$, then there is a natural exact sequence
\[
0\to HC_{2n-1+j}(\M^j)\to K_{2n+j}(\M^j)\to\Z\overset{\alpha^j_n}\to HC_{2n-2+j}(\M^j)\to K_{2n-1+j}(\M^j)\to 0
\]
If in addition $\I\subset \cL_p$ for some $p\ge 1$, then $\alpha_n^j$ is injective for all $n\ge (p+1-j)/2$.
\end{enumerate}
\end{propo}

\begin{proof}
The proof of \cite[6.4.1]{wt} shows that if $H_*$ is as in i) and $\Gamma$ is an infinite sum ring, then $H_*(\Gamma)=0$.
In particular this applies to $\Gamma=\B$. From this observation and the exact sequence
\begin{equation}\label{seqm0}
0\to M_2(\I)\to \M^0(\I)\to \B\to 0
\end{equation}
we get $H_*(\I)=H_*(\M^0(\I))$. Next, consider the map of exact sequences
\begin{equation}\label{seqm1}
\xymatrix{0\ar[r]&M_2(\I)\ar[r]\ar[d]&\M^1(\I)\ar[r]\ar[d]&\B/\I\oplus\B/\I\ar[d]\ar[r]&0\\
          0\ar[r]&M_2(\I)\ar[r] &M_2(\B)\ar[r]&M_2(\B/\I)\ar[r]&0}
\end{equation}
Comparing the long exact sequences $H_*$ associates to the top and bottom rows, we obtain that $H_*(\M^1(\I))=H_{*+1}(\I)$. This proves
i). Part ii) follows using \eqref{seqm0}, \ref{lem:infisum}, \ref{promv}, and the top row of \eqref{seqm1}.
The exact sequence of iii) follows from \ref{propo:nuiso} and \ref{propo:kopi} i). To prove that $\alpha_n^j$ is injective for
$n\ge (p+1)/2$, it suffices to consider the case $\I=\cL_p$, which follows from \cite[\S2.9]{ck}.
\end{proof}

\begin{rem} Let $\I$ be a sub-harmonic ideal, $j=0,1$, and $n\ge 1$. Assume that the map $\alpha^j_n$ of part iii) of Proposition
\ref{propo:mji} above is injective. Assume further that a group homomorphism $\phi:HC_{2n-2+j}(\M^j)\to \C$ is given such that
$\phi\circ\alpha^j_n$ is injective. Then $\phi$ induces a group homomorphism
\begin{equation}
\eta(\phi):K_{2n-1+j}(\M^j)\to \C/\alpha^j_n(\Z)\cong \C^\times.
\end{equation}
We remark that the multiplicative character of Connes and Karoubi \cite{ck} is a particular case of this construction for
$\I=\cL_{p}$ when $p=2n-1+j\ge 1$. They construct a continuous cyclic cocycle $\tau\in HC^{2n-2+j}_{\top}(\M^j)$ and show that, for the characters
of Subsection \ref{subsec:compaktop}, the composite of $\tau$ with $\psi:=S\circ ch_{2n+j}^{\dif'}$ is nonzero. Then they define
their multiplicative character as the induced map
\[
\theta:K_{2n-1+j}(\M^j)\to \C^\times.
\]
By \ref{thm:main} iii) and \eqref{compa_diags} $\psi$ is the composite of $\alpha^j_n$ with
$\iota:HC_{2n-2+j}(\M^j)\to HC^{\top}_{2n-2+j}(\M^j)$. Thus for $\phi=\tau\iota$, we have $\theta=\eta(\phi)$.
\end{rem}

\section{Karoubi's conjecture}\label{sec:karconj}

Several conjectures circulated under the name \textit{Karoubi's Conjecture} in the literature. One conjecture concerning the
algebraic $K$-theory of stable $C^*$-algebras (formulated by Higson in \cite[6.5]{hig})
was resolved by A. Suslin and M. Wodzicki in \cite{qs}. The proof is based on their
result characterizing excision in algebraic $K$-theory in terms of $H$-unitality and on the homotopy invariance theorem
for $C^*$-algebras \cite{hig}, which is in turn based on work of G. Kasparov, N. Higson and J. Cuntz.
Another conjecture, made by Karoubi in \cite{karcomp}, was about the algebraic $K$-theory of unital Banach algebras, stabilized with
the algebra of compact operators, using the projective tensor product. M. Wodzicki is credited with the solution of a strong version of
the latter conjecture. His result is stated without proof in \cite[Theorem 2 (a)]{wod}. He has lectured on his proof;
however the latter has not appeared in the literature.
In \ref{thm:karconj} we state and give a complete proof of this result. Our argument uses both results of Wodzicki whose proof has been published, and
results of this paper which are independent of his work.

\subsection{Homological preliminaries}\label{subsec:homprep}
Throughout this section, $k$ is a field of characteristic zero.
Let $A$ be a $k$-algebra, $\tilde{A}_k$ its unitalization.
Recall from \cite[\S2]{wodex} that the {\it bar homology} of $A$ is
\[
H^{\rm bar}_*(A/k):=\tor^{\tilde{A}_k}_{*}(k,A).
\]
We write $\tilde{A}$ and $H_*^{\rm bar}(A)$ for $\tilde{A}_\Q$ and $H^{\rm bar}_*(A/\Q)$.

\begin{lem} Let $A$ be a $k$-algebra. Then $H^{\rm bar}_*(A/k)=H^{\rm bar}_*(A)$.
\end{lem}
\begin{proof}
Note that
\begin{gather*}
\tilde{A}_k\otimes_{\tilde{A}}N=k\otimes N\qquad (N\in \mod_{\tilde{A}})\\
k\otimes_{\tilde{A}_k}M=M/AM=\Q\otimes_{\tilde{A}}M\qquad (M\in\mod_{\tilde{A}_k}).
\end{gather*}
It follows from this that if $L\weq A$ is an $\tilde{A}_k$-free resolution, then it is also an $\tilde{A}$-flat resolution, and
$H_*^{\rm bar}(A/k)=H_*(L/AL)=H_*^{\rm bar}(A)$.
\end{proof}

\begin{rem}
Recall the algebra $A$ is called {\it $H$-unital} if $H^{bar}_*(A/k)=0$. By the previous lemma, this definition is independent of the field
$k\supset \Q$.
\end{rem}

The following lemma is a particular case of a spectral sequence for Hochschild homology which has been observed by various people,
including C. Weibel and J.A. and J.J. Guccione (cf. personal communication with the first author).

\begin{lem} Let $R$ be a unital $k$-algebra, and $M$ an $R\otimes_kR^{op}$-module. Write
$HH^k_*(?)$ for Hochschild homology of $k$-algebras. Consider the graded algebra $\Omega^*_{k}$ of (absolute, commutative)
K\"ahler differential forms of $k$.
There is a spectral sequence
\[
E^2_{pq}=HH^k_p(R,M)\otimes\Omega^q_{k}\Rightarrow HH_{p+q}(R,M).
\]
\end{lem}
\begin{proof}
The functor $\mod_{R\otimes R^{op}}\to \mod_k$, $M\mapsto HH_0(M)=M/[M,R]$ factors as
$HH_0(k,?):\mod_{R\otimes R^{op}} \to \mod_{R\otimes_k R^{op}}$
followed by $HH_0^k(R,?):\mod_{R\otimes_k R^{op}}\to \mod_k$. A computation shows that the $E^2$-term of the Grothendieck
spectral sequence is as stated in the lemma.
\end{proof}

\begin{coro}\label{coro:0k=0} Let $A$ be a $k$-algebra. Then $HH^k_*(A)=0$ if and only if $HH_*(A)=0$.
\end{coro}
\begin{proof} The direction $\Rightarrow$ is clear. The converse follows from an inductive argument, using the spectral sequence
above.
\end{proof}

\begin{lem}\label{lem:kuenneth} Let $A$ and $B$ be $H$-unital $k$-algebras. Then $HH_*^k(A\otimes_k B)=HH^k_*(A)\otimes_kHH^k_*(B)$.
\end{lem}
\begin{proof} If $A$ and $B$ are unital, this is well-known. The general case follows from this and excision.
\end{proof}

\begin{coro}\label{coro:a0=0} If $A$ and $B$ are as in the previous lemma, and $HH_p^{k}(B)=0$ for $0\le p\le n$, then $HH_p(A\otimes_kB)=0$
for $0\le p\le n$.
\end{coro}

\subsection{Homology of operator ideals}

The following theorem is an equivalent formulation of \cite[Thm. 4]{wod} (the original statement of {\it loc. cit.} is recalled
in \ref{rem:actual} below). Since no proof of \cite[Thm. 4]{wod} has appeared in print, we give a short argument, referring
only to results with a published proof.

\begin{thm}\label{thm:jhunit}\cite[Thm. 4]{wod}
Let $\J$ be an operator ideal. If $\J=\J^2$, then $\J$ is $H$-unital.
\end{thm}
\begin{proof}
In \cite[Sec. 3]{qs} it was shown that an algebra $\J$ which satisfies the triple factorization property is $H$-unital. Recall that an algebra $\J$
is said to satisfy the (right) triple factorization property, if the following condition is satisfied:

For any finite collection $a_1,\dots,a_n \in \J$, there exist $b_1,\dots,b_n,c,d \in \J$, such that the following two conditions
hold.
\begin{enumerate}
\item[(i)] $a_i=b_icd,\quad \forall i \in \{1,\dots,n\}$.
\item[(ii)] The left annihilator of $c$ coincides with the left annihilator of $cd$.
\end{enumerate}
We proceed by showing that the equality $\J^2=\J$ implies that $\J$ has the triple factorization property.

Let $a_1, \dots a_n \in \J$ be any finite collection and $a_i= u_i b'_i$ the polar decomposition.
Consider the element
\[c =\left(b'_1 + \dots + b'_n \right)^{\frac18}. \]
Because
\[
b'_i=({{b'}_i}^{\frac12})^*\cdot {b'}_i^{\frac12}\le c^8
\]
there exists a bounded operator $c_i$ such that ${b'_i}^{\nicefrac 12}  = c_i c^2$ (\cite[I.4.6]{dav}).
Now, put $b_i = u_i {b'_i}^{\nicefrac12} c_i$ and compute $a_i = b_i c^2$. Note
that the left annihilator of $c^2$ conincides with that of $c$ since $c$ is positive. Note further
that $b_i, c \in \J$ since by assumption $\J^{2^{-n}}=(\J^{2^n})^{2^{-n}} =\J$.
\end{proof}

\begin{rem}\label{rem:actual} The precise claim of \cite[Theorem 4]{wod} is that if $\J$ is an operator ideal with $\J=\J^2$, and $A$
an $H$-unital $\C$-algebra, then $A\otimes_\C\J$
satisfies excision in algebraic $K$-theory. This follows from \ref{thm:jhunit}, combining the following results.
In \cite[Thm. 7.10]{qs}, Suslin and Wodzicki showed that the category of $H$-unital algebras over a fixed ground ring $k$ is closed
under tensor products. Furthermore, they proved
in \cite[Thm. B]{qs}, that a $\Q$-algebra satisfies excision in algebraic $K$-theory if and only if it is $H$-unital.
\end{rem}

\begin{lem}\label{lem:hhj=0}
Let $\J$ be an operator ideal that is either sub-harmonic or multiplicative. Then the following conditions are equivalent:
\begin{enumerate}
\item[(i)] \quad $\J=[\J,\J]$.
\item[(ii)] \quad $HH_*(\J)=0$.
\end{enumerate}
\end{lem}
\begin{proof}
The only nontrivial implication is (i) $\Rightarrow$ (ii). Let us assume first that $\J$ is multiplicative.
A choice of an isomorphism $H \otimes_2H \cong H \oplus H$, which identifies $H \otimes e_1$ with $H \oplus 0$ induces
a commutative diagram
\[
\xymatrix{\J\ar[dr]_{1_{\J} \otimes e_{11}}\ar[rr]^{1_{\J} \oplus 0}&&M_2(\J)\\
                            &\J\otimes \J\ar[ur]&}
\]
Condition (i) implies that $\J = \J^2$. Hence, by \ref{thm:jhunit}, $\J$ is $H$-unital
as a $\Q$-algebra. Thus the horizontal map induces an isomorphism
in $HH_*$. Using Lemma \ref{lem:kuenneth} and induction, we get that $HH_n(\J)=0$ for all $n$.

If $\J$ is sub-harmonic, then $\cL_p$ multiplies $\J_\infty$ for some $p>0$. Since $\J^2 = \J$, we actually get that the multiplicative ideal
$(\cL_1)_{\infty}$ multiplies $\J$. Since $[(\cL_1)_{\infty},(\cL_1)_{\infty}] = (\cL_1)_{\infty}$ (as follows from \cite[Thm. 2]{PT}),
the above result applies to $(\cL_1)_{\infty}$. The proof can be finished by arguing as above, but using the following diagram instead:
\[
\xymatrix{\J\ar[dr]_{1_{\J} \otimes e_{11}}\ar[rr]^{1_{\J} \oplus 0}&&M_2(\J)\\
                            &\J\otimes (\cL_1)_{\infty}.\ar[ur]&}
\]
\end{proof}

\begin{rem}\label{rem:conj1}
Conjecturally, the preceding result is true for all operator ideals.
\end{rem}

\begin{thm}\label{thm:alg_carcoj}
Let $A$ be an H-unital $\C$-algebra and let $\J$ be a sub-harmonic ideal satisfying $[\J,\J]=\J$.
\begin{enumerate}
\item[(i)] The cyclic homology groups $HC_*(A \otimes_\C \J)$ vanish.
\item[(ii)] The algebra $A \otimes_\C \J$ is $K$-regular.
\item[(iii)] The graded abelian group $K_*(A \otimes_\C \J)$ is $2$-periodic.
\end{enumerate}
\end{thm}
\begin{proof}
Part i) follows from \ref{coro:0k=0}, \ref{coro:a0=0}, \ref{thm:jhunit}, and \ref{lem:hhj=0}. Parts ii) and iii) follow from i) and \eqref{absoseq}.
\end{proof}

\begin{rem}
Note that the operator ideals $(\cL_1)_\infty$ and $\cL_\infty$ are examples of operator ideals satisfying the assumption of the preceding
corollary. Every other such ideal has to lie between these. The particular case of part iii) of \ref{thm:alg_carcoj} when
$\J=\cL_\infty$ is \cite[Theorem 2 (c)]{wod}.
\end{rem}
\subsection{Karoubi's conjecture}

The following result relating functional analytic properties to purely algebraic properties has been obtained by Wodzicki in \cite[Thm. $8.1$]{wodex}.
\begin{thm}[Wodzicki] \label{thm:wodzh}
Let $A$ be a Fr\'echet algebra, whose topology is generated by a countable family of sub-multiplicative seminorms
$\{\rho_n: A \to \R_+, n \in \N\}$. Assume further that there exists an approximate right or left unit in $A$,
which is totally bounded with respect to the family of seminorms. Under these circumstances, the algebra $A$ is $H$-unital.
\end{thm}
\begin{exa} The operator ideal $\cL_{\infty}$ of all compact operators is $H$-unital since it has a bounded approximate unit and thus satisfies the
assumptions of Theorem \ref{thm:wodzh}. Similarly, if $L \in \locass$ satisfies the hypothesis of Theorem \ref{thm:wodzh}, then so does
$L \hotimes \cL_\infty$.
\end{exa}
\sn
Next we prove Wodzicki's theorem that Karoubi's conjecture holds for stable Fr\'echet algebras satisfying the condition
of the theorem above. The referee informs us that Wodzicki's original argument went along the line that we present here.
\sn
\begin{thm}[Karoubi's Conjecture]\label{thm:karconj}
Let $L$ be a Fr\'echet algebra, satisfying the hypothesis of Theorem \ref{thm:wodzh}. There is a natural isomorphism:
\[ K_n(L \hotimes \cL_\infty) \stackrel{\sim}{\to} K_n^{\top}(L \hotimes \cL_\infty), \quad \forall n \in \Z.\]
\end{thm}
\begin{proof}
We claim that $HH_*(L \hotimes \cL_\infty)=0$. This implies that $HC_*(L \hotimes \cL_\infty)=0$ and the proof is immediate from this and
\ref{thm:wodzh}, using the exact sequence of Theorem \ref{thm:exact}. To prove the claim, note
that, by Corollary \ref{coro:0k=0}, the latter is equivalent to the vanishing of
$HH^\C_*(L \hotimes \cL_\infty)$ . Next note that the map $HH^{\C}_*(L \hotimes \cL_\infty) \to HH^{\C}_*(L \hotimes M_2(\cL_\infty))$
induces the identity and factors through
$HH^\C_*(L \hotimes \cL_\infty \otimes_\C \cL_\infty) \cong HH^\C_*(L \hotimes \cL_{\infty}) \otimes_\C HH^\C_*(\cL_{\infty})$.
Here, we used $H$-unitality and Lemma \ref{lem:kuenneth}. The proof is finished by recalling that $HH^\C_*(\cL_\infty)=0$
by Lemma \ref{lem:hhj=0} and Corollary \ref{coro:0k=0}.
\end{proof}

\begin{rem}\label{rem:kartriv}
Note that a direct proof of Karoubi's conjecture does not require the whole machinery that we have developed. The key observation is
the strength of Theorem \ref{thm:wodzh}. We sketch an alternative argument. Theorem \ref{thm:wodzh} implies
that for $L$ as in \ref{thm:wodzh}, the functor from $C^*$-algebras to abelian groups
\[
A \mapsto K_*(L \hotimes (A \otimes_{\rm min} \cL_\infty))
\]
is excisive and $\cL_{\infty}$-stable. Using the homotopy invariance theorem for $C^*$-algebras (see \cite{hig})
we obtain that such a functor is homotopy invariant. It follows from this that the restriction of
$L\mapsto K_*(L\hotimes\cL_\infty)$ to the full subcategory $\mathfrak{C}\subset\locass$ of those locally convex algebras which satisfy the hypothesis of
\ref{thm:wodzh} is diffeotopy invariant. Next note that, because $C^\infty(\Delta^p)$ is unital ($p\ge 0$), the functor
 $?\hotimes C^\infty(\Delta^p)$ maps $\mathfrak{C}$ to itself. Thus for $L\in\mathfrak{C}$,
\[
K_*(L\hotimes\cL_\infty)=K_*(L\hotimes C^\infty(\Delta^p)\hotimes\cL_\infty).
\]
It follows that
\[
K_*(L\hotimes\cL_\infty)=KD_*(L\hotimes\cL_\infty)=K^{\top}_*(L\hotimes\cL_\infty).
\]
This concludes the proof. We remark further that, if we just want to prove the result for Banach algebras with bounded approximate units,
the proof simplifies, as we can skip the step from continous to differentiable homotopy invariance.

The key fact used in the short argument above is that all $C^*$-algebras have bounded approximate units. However, note that the proof
presented in \ref{thm:karconj} is more conceptual, since
it shows directly that the obstruction to having an isomorphism vanishes, and produces interesting results along the way, which leave room for
a generalization of Karoubi's Conjecture beyond $\cL_\infty$. We conjecture that, in an appropriate setting of bornological algebras, Karoubi's Conjecture holds for
$(\cL_1)_{\infty}$ in place of $\cL_\infty$.
\end{rem}

\end{document}